\documentclass{article}
\usepackage{amsthm}
\usepackage[T1]{fontenc}
\usepackage[frenchb, english]{babel}
\usepackage[hmargin=4cm,vmargin=4cm]{geometry}
\usepackage{amssymb}
\usepackage{amsmath}

\def\aro{{A}^{\!\!\!\raise5pt\hbox{$\scriptstyle \circ$}}}
\def\aroi{{A}^{\!\!\!\raise4pt\hbox{$\scriptscriptstyle \circ$}}}
\def\cro{\smash{{C}^{\!\!\!\raise5pt\hbox{$\scriptstyle \circ$}}}}
\def\croi{\smash{{C}^{\!\!\!\raise4pt\hbox{$\scriptscriptstyle \circ$}}}}

\begin{document}

\newtheorem{Th}[subsubsection]{Théorème}
\newtheorem{Pro}[subsubsection]{Proposition}
\newtheorem{De}[subsubsection]{Définition}
\newtheorem{Prt}[subsubsection]{Propriété}
\newtheorem{Prts}[subsubsection]{Propriétés}
\newtheorem{Le}[subsubsection]{Lemme}
\newtheorem{Hyp}[subsubsection]{Hypothèse}
\newtheorem{Cor}[subsubsection]{Corollaire}



\title{Cramér's theorem\\
in measurable locally convex spaces}
\date{March 17, 2011}
\author{Pierre Petit\\ \\
Universit\'e Paris Sud}
\maketitle

\selectlanguage{english}

\begin{abstract}
We give a general setting for Cramér's large deviations theorem for the empirical means $\overline{X}_n$ of a sequence $(X_n)_{n \geqslant 1}$ of i.i.d. random vectors, which contains Cramér's theorem in a Banach space and Sanov's theorem. We define the notion of \emph{measurable locally convex space} in which the variables $\overline{X}_n$ are indeed measurable. We obtain the standard weak large deviations principle for the sequence $(\overline{X}_n)_{n \geqslant 1}$ and, having a close look at the notion of convex tension, we prove that the upper bound holds for convex sets. Hence we show that the identification between the entropy and the opposite of the Fenchel-Legendre transform of the pressure, $s = -p^*$, holds in any measurable locally convex space. The proof is based on convex duality and monotone convergence and does not resort to the law of large numbers or any other limit theorem. We also show that, if $\mu$ is the law of $X_1$, then $\overline{\textup{dom}(s)} = \textup{cosupp}(\mu)$.
\end{abstract}

\selectlanguage{french}

\begin{abstract}
Nous établissons un cadre général pour le théorème de Cramér sur les grandes déviations des moyennes empiriques $\overline{X}_n$ d'une suite $(X_n)_{n \geqslant 1}$ de vecteurs aléatoires i.i.d., cadre qui contient le théorème de Cramér dans les espaces de Banach séparables et le théorème de Sanov. Nous introduisons la notion d'\emph{espace vectoriel localement convexe mesurable} dans lequel les variables $\overline{X}_n$ sont bien mesurables. On obtient le principe de grandes déviations faible classique pour la suite $(\overline{X}_n)_{n \geqslant 1}$ ainsi que, en examinant finement la notion de convexe-tension, la borne supérieure pour les convexes. Ainsi, on montre que l'identification entre l'entropie et l'opposée de la transformée de Fenchel-Legendre de la pression, $s = -p^*$, est vraie dans tout espace vectoriel localement convexe mesurable. La preuve repose sur la dualité convexe et la convergence monotone, et ne fait appel ni à la loi des grands nombres ni à un autre théorème-limite. On montre aussi que, notant $\mu$ la loi de $X_1$, $\overline{\textup{dom}(s)} = \textup{cosupp}(\mu)$.
\end{abstract}

\medskip

\textbf{Remarque :} La prochaine version de ce texte sera en langue anglaise.

\section{Introduction}

Etant donné une suite $(X_n)_{n \geqslant 1}$ de variables aléatoires i.i.d. à valeurs dans un espace vectoriel, le théorème de Cramér dit que la suite des mesures empiriques $\overline{X}_n$ vérifie un principe de grandes déviations (PGD) faible. La version originelle de ce théorème, dans $\mathbb{R}$, une fois sa forme générale démontrée par Bahadur \cite{Bah71}, a été rapidement généralisé en dimension supérieure. Le cadre des espaces vectoriels localement convexes a été introduit par \cite{BaZ79} et semble le cadre le plus naturel pour la théorie de Cramér. L'objet du présent texte est, au départ, une exploration de la question de l'identification entre l'entropie $s$ de la suite $(X_n)_{n \geqslant 1}$ et l'opposée de la transformée de Fenchel-Legendre de la pression $p$, autrement dit de l'égalité
$$
s = -p^*
$$
et une recherche de contre-exemples. Cela a conduit à définir un cadre, le plus général possible, ou la question a un sens : il s'agissait de résoudre les problèmes de mesurabilité des moyennes empiriques $\overline{X}_n$, d'une part, et des formes linéaires $\lambda$ sur l'espace, d'autre part. Nous ne travaillons plus avec une tribu borélienne. D'ailleurs, dans le cas du théorème de Sanov, la tribu cylindrique n'est pas la tribu borélienne de la topologie faible (ce n'est vrai que si l'espace sous-jacent est polonais). Cela nous amène à introduire le concept d'\emph{espace vectoriel localement convexe mesurable} (\emph{e.v.l.c.m.}) où la séparabilité assure la mesurabilité voulue. La tribu considérée est la \og{}tribu des convexes\fg{} (analogue de la tribu des boules, dans un espace métrique, et qui peut être plus petite que la tribu borélienne ; \emph{cf.} \cite[exemple 1.4.]{Bil68}) qui semble plus naturelle dans la théorie de Cramér : si $\mathcal{F}$ est une tribu quelconque, notre définition de l'entropie ne dépend que des valeurs de la mesure $\mu$, loi de $X_1$, sur la tribu engendrée par les convexes mesurables.

\medskip

Nous introduisons la notion de probabilité $\mu$ \emph{portée} par un ensemble $D$ : cela revient à dire que $\mu$ est la loi d'une variable à valeurs dans $D$. Dire qu'une mesure est portée par $D$ est plus fort que de dire que son support $\textup{supp}(\mu)$ est inclus dans $D$ (il existe des mesures de probabilité de support vide, et pourtant non portées par le vide). Au départ, il s'agit simplement de clarifier un point de \cite[chapitre 24]{Cer07} : le PGD faible pour une suite de mesures portées par un même ensemble relativement compact n'assure pas le PGD.
Au passage, le PGD faible en tribu cylindrique est vide si l'espace est trop gros : par exemple, dans $\mathbb{R}^\mathbb{R}$, il n'y a pas de mesurable relativement compact car tous les mesurables non vides sont non bornés. Pourtant, on s'attend à avoir le PGD si la mesure $\mu$ est portée par une partie $D$ relativement compacte (cas du théorème de Sanov). C'est pour ce genre d'exemples, où l'argument de sous-aditivité continue à fonctionner, que nous avons rajouté les notions de probabilité portée par un ensemble et de borne supérieure \textsf{(BS$_{\flat , D}$)} pour les ensembles dont l'intersection avec $D$ est relativement compacte. Ces considérations s'avèrent ensuite utiles pour passer de l'égalité $s = -p^*$ dans l'espace séparable et complet $\mathcal{C}([0, 1] ; \mathbb{R})$ à la même égalité dans tout sous-espace. Or, les théorèmes de Banach-Mazur et de Mazur-Ulam permettent de voir que tout espace vectoriel normé séparable est un sous-espace de $\mathcal{C}([0, 1] ; \mathbb{R})$. On en déduit l'égalité $s = -p^*$ dans tout e.v.l.c.m., en toute généralité.

\medskip

Enfin, nous empruntons à \cite{Mor67} l'idée de travailler avec des fonctions à valeurs dans $[-\infty , +\infty]$. Dans ce cadre, une fonction convexe est automatiquement propre, ou bien constante et égale à $\pm \infty$. Cela simplifie la terminologie et évite des vérifications de finitude (sur les transformées de Fenchel-Legendre) superflues.

\medskip

Le déroulement de la preuve est le suivant (et analogue à celui de \cite{CeP10} : le lemme sous-additif et le \emph{principle of the largest term} permettent d'établir le PGD faible. Nous montrons même, plus généralement, un lemme de Varadhan compact \ref{varadcomp}. Puis, la \emph{convexe-tension} nous permet de passer de la borne supérieure pour les compacts \textsf{(BS$_\flat$)} à la borne supérieure pour les convexes \textsf{(BS$_c$)}. En cela, nous simplifions la preuve de \cite{BaZ79} et nous nous passons de l'hypothèse additionnelle 6.1.2 (b) de \cite{DeZ93}, introduite pour rectifier la preuve de \cite[appendice, proposition 1]{BaZ79}. La démonstration de $p = (-s)^*$, que nous donnons avant, n'est qu'un cas particulier du lemme de Varadhan convexe \ref{varadconv}. Même si cela n'apparaît pas, car nous donnons une preuve directe de l'égalité $p = (-s)^*$, il s'avère que la borne supérieure pour les demi-espaces \textsf{(BS$_h$)} entraîne l'égalité $p = (-s)^*$ (cela apparaît dans la preuve de \cite{BaZ79}, puis de \cite{Aze80}, et la remarque est faite dans \cite[proposition 16.2.]{Cer07}). Ensuite, l'égalité $s = -p^*$ découle de la propriété d'inversion de la transformée de Fenchel-Legendre.

\medskip

Nous nous intéressons enfin à d'autres questions apparaissant naturellement dans la théorie de Cramér, la première d'entre elles étant la description du domaine $\textup{dom}(s)$ de l'entropie. L'hypothèse de convexe-tension nous suffit à montrer que
$$
\overline{\textup{dom}(s)} = \textup{cosupp}(\mu)
$$
Ainsi, nous nous passons de l'hypothèse de régularité de \cite{BaZ79} et approfondissons les résultats de \cite{Cer07} (notamment lemme 9.7.\footnote{Remarquons que la démonstration de \cite[lemme 9.8.]{Cer07} suppose la convexe-régularité.}). Nous nous intéressons également à l'existence de la limite
$$
\lim_{n \to \infty} \frac{1}{n} \log \mu_n(C)
$$
pour les convexes mesurables $C$ (\emph{cf.} \cite{Sla88} pour un bon résumé des résultats pour des ensembles $C$ non nécessairement convexes). Dans les suppléments techniques, nous évoquons la question de la séparation : elle n'est pas nécessaire pour la démonstration, elle simplifie les termes employés. En fait, la borne supérieure faible est vraie pour des ensembles vérifiant une propriété de sous-recouvrement fini par des convexes ouverts : de même que la convexe-régularité s'avère plus pertinente que la régularité, de même cette propriété de sous-recouvrement fini par des convexes ouverts est plus pertinente que la propriété de Borel-Lebesgue. Ce genre de question est aussi derrière les problèmes qui nous ont amenés à introduire l'ensemble $D$. Il est sans doute possible de faire une théorie de Cramér sans parler de compacts, mais en identifiant des ensembles plus naturels, vérifiant une propriété convexe de Borel-Lebesgue.




\section{Cadre}

\subsection{Entropie et pression}

Soient $X$ un espace vectoriel réel, $\mathcal{F}$ une tribu sur $X$ et $\tau$ une topologie séparée (au sens de Hausdorff) sur $X$. Soit $(X_n)_{n \geqslant 1}$ une suite de variables aléatoires i.i.d. à valeurs dans $X$. On définit ses moyennes empiriques par
$$
\overline{X}_n := \frac{1}{n} \sum_{k=1}^n X_k
$$
L'\emph{entropie} de $(X_n)_{n \geqslant 1}$ est la fonction $s : X \rightarrow [-\infty, 0]$ définie par
$$
\forall x \in X \qquad s(x) := \inf_{\substack{A \in \mathcal{F} \\ \aroi \ni x}} \liminf_{n \to \infty} \frac{1}{n} \log \mathbb{P}\big( \overline{X}_n \in A \big)
$$
Par construction, l'entropie $s$ est la plus grande fonction vérifiant la borne inférieure :

\medskip

\textsf{(BI)} pour tout $A \in \mathcal{F}$,
$$
\liminf_{n \to \infty} \frac{1}{n} \log \mathbb{P}\big( \overline{X}_n \in A \big) \geqslant \sup_{\aroi} s
$$

On dit que $(X_n)_{n \geqslant 1}$ vérifie un \emph{principe de grandes déviations} (PGD) si la borne supérieure suivante est vérifiée :

\medskip

\textsf{(BS)} pour tout $A \in \mathcal{F}$,
$$
\limsup_{n \to \infty} \frac{1}{n} \log \mathbb{P}\big( \overline{X}_n \in A \big) \leqslant \sup_{\smash{\overline{A}}} s
$$

De manière générale, \textsf{(BS)} n'est pas vérifiée pour tous les mesurables. Si $\mathcal{P}$ désigne un ensemble de parties mesurables de $X$, on définit la version restreinte de la borne supérieure suivante :

\medskip

\textsf{(BS$_\mathcal{P}$)} pour tout $A \in \mathcal{P}$,
$$
\limsup_{n \to \infty} \frac{1}{n} \log \mathbb{P}\big( \overline{X}_n \in A \big) \leqslant \sup_{\smash{\overline{A}}} s
$$

En particulier, si $D$ est une partie de $X$, on notera :

\medskip

\textsf{(BS$_{\flat , D}$)} pour tout $K \in \mathcal{F}$ tel que $K \cap D$ soit relativement compact,
$$
\limsup_{n \to \infty} \frac{1}{n} \log \mathbb{P}\big( \overline{X}_n \in K \big) \leqslant \sup_{\smash{\overline{K}}} s
$$

\textsf{(BS$_{c, D}$)} pour tout $C \in \mathcal{F}$ tel que $C \cap D$ soit convexe,
$$
\limsup_{n \to \infty} \frac{1}{n} \log \mathbb{P}\big( \overline{X}_n \in C \big) \leqslant \sup_{\smash{\overline{C}}} s
$$

Si $(X_n)_{n \geqslant 1}$ vérifie \textsf{(BS$_{\flat , X}$)}, on dit que $(X_n)_{n \geqslant 1}$ vérifie un \emph{principe de grandes déviations faible} (PGD faible). L'objectif principal de ce texte est de donner un cadre général simple pour avoir un PGD faible et une expression simple de l'entropie en fonction de la loi de $X_1$. Explicitons ce dernier point. Notons $X^*$ le dual topologique de $X$. La \emph{pression} de $(X_n)_{n \geqslant 1}$ est l'application $p : X^* \rightarrow [-\infty, +\infty]$ définie par
$$
\forall \lambda \in X^* \qquad p(\lambda) := \log \mathbb{E} \big( e^{\langle \lambda | X_1 \rangle } \big)
$$
Dans notre cadre, l'entropie est l'opposée de la fonction convexe-conjuguée (aussi appelée transformée de Fenchel-Legendre) de la pression et a donc pour expression
$$
-p^*(x) := - \sup_{\lambda \in X^*}  \big( \langle \lambda | x \rangle - p(\lambda) \big)
$$
On notera que les définitions mêmes de $s$ et de $p$ posent un problème de mesurabilité. Le cadre général que nous proposons ci-après résout ce problème. Nous ne savons pas si ce cadre contient le cadre de \cite{BaZ79} et de \cite{Cer07} (qui prennent la mesurabilité comme hypothèse), en tout cas il est explicite et plus général que leurs applications (aux espaces de Banach séparables, aux espaces polonais, ou au théorème de Sanov en $\tau$-topologie et tribu cylindrique).

\subsection{Espaces vectoriels localement convexes mesurables}

Soit $X$ un espace vectoriel réel. On se propose de munir $X$ d'une topologie et d'une tribu de sorte que l'entropie et la pression soient définies pour toute suite de variables i.i.d. donnée sur $X$. Si $C$ est un convexe de $X$ contenant $0$, sa jauge (ou fonctionnelle de Minkowski) est l'application $M_C : X \rightarrow [0, +\infty]$ définie par
$$
\forall x \in X \qquad M_C(x) = \inf \{ t \geqslant 0 ; x \in tC \}
$$
La proposition suivante relie les propriétés topologiques d'un ouvert convexe à ses propriétés algébriques intrinsèques.
\begin{Prt}\label{convj}
Soient $(X, \tau)$ un espace vectoriel topologique réel et $C$ un voisinage convexe de $0$. Alors $M_C$ est finie partout et
$$
\cro = \{ x \in X ; M_C(x) < 1 \} \qquad \textrm{et} \qquad \overline{C} = \{ x \in X ; M_C(x) \leqslant 1 \}
$$
\end{Prt}
Pour la preuve, on renvoie par exemple à \cite[II.21]{BouEVT}. Pour les deux égalités, on peut utiliser la continuité de $M_C$ pour montrer que $M_{\croi} = M_C = M_{\smash{\overline{C}}}$. Cette proposition justifie la définition suivante : un convexe $C$ de $X$ est un \emph{convexe interne} si sa jauge $M_C$ est finie partout et si
$$
C = \{ x \in X ; M_C(x) < 1 \}
$$
La notion de convexe interne est centrale dans cet article : elle fait le lien entre la structure topologique, l'outil \og{}jauge\fg{} qui permet de mener efficacement les calculs, et surtout le lemme sous-additif qui permet de définir l'entropie. Si $C$ est un convexe interne, $\{ tC \, ; \, t \in \mathbb{R}^* \}$ est un système de voisinages de $0$ pour une topologie localement convexe (non nécessairement séparée) $\tau_C$ sur $X$ (cf. \cite[II.25]{BouEVT}) : on dira que $\tau_C$ est \emph{la topologie localement convexe engendrée par $C$}. On dit qu'un convexe $C$ est \emph{symétrique} si $C = -C$. On se donne maintenant une famille $\mathcal{C}_0$ de convexes de $X$ vérifiant les trois axiomes suivants :

\medskip

\textup{\textsf{(EVLCM$_1$)}} pour tout $C \in \mathcal{C}_0$, $C$ est un convexe interne symétrique\footnote{L'hypothèse symétrique n'est pas fondamentale, mais on peut toujours s'y ramener, quitte à remplacer les $C \in \mathcal{C}_0$ par $C \cap (-C)$.} et la topologie localement convexe $\tau_C$ engendrée par $C$ est séparable ;

\medskip

\textup{\textsf{(EVLCM$_2$)}} $\mathcal{C}_0$ est stable par intersection finie et par dilatation de rapport non nul ;

\medskip

\textup{\textsf{(EVLCM$_3$)}} L'intersection des éléments de $\mathcal{C}_0$ est réduite à $\{ 0 \}$.

\medskip

Alors $\mathcal{C}_0$ est un système fondamental de voisinages de $0$ pour une topologie localement convexe séparée\footnote{C'est l'axiome \textup{\textsf{(EVLCM$_3$)}} qui assure la séparation. L'hypothèse n'est pas nécessaire (cf. Suppléments techniques) mais nous la faisons pour deux raisons. La première d'ordre culturel : pour certains auteurs, la séparation fait partie des axiomes d'une topologie d'espace vectoriel topologique. La seconde d'ordre pratique : nous jugeons plus clair de manipuler des compacts que des quasi-compacts.} $\tau$ sur $X$ (\emph{cf.} \cite[II.25]{BouEVT}). Pour ce qui est de la tribu, notons, pour tout $x \in X$,
$$
\mathcal{C}_x = \{ x + C \, ; \, C \in \mathcal{C}_0 \}
$$
un système fondamental de voisinages convexes ouverts de $x$. On définit
$$
\mathcal{F} = \sigma\left( \bigcup_{x \in X} \mathcal{C}_x \right)
$$
On dira alors que le quadruplet $(X, \mathcal{C}_0, \mathcal{F}, \tau)$ est un \emph{espace vectoriel localement convexe mesurable} (\emph{e.v.l.c.m.}). On notera que $\mathcal{C}_0$ est un système fondamental de voisinages de $0$ convexes, symétriques, mesurables et ouverts ; et ils vérifient donc la propriété ci-dessus. On dira que $(X, \mathcal{C}_0, \mathcal{F}, \tau)$ est \emph{l'e.v.l.c.m. associé} à $\mathcal{C}_0$. Deux remarques pour finir ce paragraphe. La première : la topologie $\tau$ n'est pas nécessairement séparable : par exemple, $X = \mathbb{R}^{2^\mathbb{R}}$ muni de la topologie produit n'est pas séparable et peut être obtenu avec la construction ci-dessus. La seconde : en pratique, on se donne une famille $\widetilde{\mathcal{C}}_0$ vérifiant \textsf{(EVLCM$_1$)} et \textsf{(EVLCM$_3$)} et on construit $\mathcal{C}_0$ en ajoutant à $\tilde{\mathcal{C}}_0$ les dilatés de rapports non nuls des éléments de $\widetilde{\mathcal{C}}_0$, puis les intersections finies des ensembles obtenus. On vérifie alors que $\mathcal{C}_0$ est bien stable par dilatation de rapport non nul et que tout élément de $\mathcal{C}_0$ est un convexe interne symétrique engendrant une topologie localement convexe séparable.

\subsection{Suites de Cramér et probabilités portées par une partie}

Soient $(X, \mathcal{C}_0, \mathcal{F}, \tau)$ un e.v.l.c.m. et $\mu$ une mesure de probabilité sur $\mathcal{F}$. Soit $(X_n)_{n \geqslant 1}$ une suite de variables indépendantes et de même loi $\mu$\footnote{L'existence d'une telle suite est assurée, en toute généralité, par un résultat énoncé par \L{}omnicki et Ulam \cite{LoU34} dont la première démonstration est due à von Neumann \cite{vNe35}. Il n'y a donc pas besoin du théorème d'extension de Kolmogorov dans le cas précis des mesures produits. On renvoie à la bibliographie de \cite{SpJ48} pour plus de détails.}. Pour tout $n \geqslant 1$, on note $\mu_n$ la mesure image de $\mu^{\otimes n}$ par l'application mesurable (cf. Questions de mesurabilité)
$$
(x_1, \ldots , x_n) \in X^n \mapsto \frac{1}{n} \sum_{k=1}^n x_k \in X
$$
autrement dit la loi de la moyenne empirique
$$
\overline{X}_n = \frac{1}{n} \sum_{k=1}^n X_k
$$
On dira que $(\mu_n)_{n \geqslant 1}$ est la \emph{suite de Cramér associée à $\mu$}. Si $D$ est une partie de $X$, on dira que $\mu$ est \emph{portée par $D$} si, pour tout $(A, B) \in \mathcal{F}^2$,
$$
A \cap D = B \cap D \Rightarrow \mu(A) = \mu(B)
$$
Dans ce cas, $\mu$ est la loi d'une variable aléatoire $X_1$ à valeurs dans $D$ : il suffit de considérer l'inclusion
$$
X_1 : (D, \mathcal{F}|_D, \mu|_D) \hookrightarrow (X, \mathcal{F})
$$
où $\mathcal{F}|_D = \{ A \cap D \, ; \, A \in \mathcal{F} \}$ et $\mu|_D(A \cap D) = \mu(A)$ pour tout $A \in \mathcal{F}$ : le fait que $\mu$ soit portée par $D$ justifie la définition de la probabilité $\mu|_D$, appelée \emph{probabilité trace de $\mu$ sur $D$}. Réciproquement, si $\mu_D$ est une probabilité sur $\mathcal{F}|_D$, on peut définir une mesure $\mu$ sur $\mathcal{F}$ par
$$
\forall A \in \mathcal{F} \qquad \mu(A) = \mu_D(A \cap D)
$$
Alors, $\mu$ est une probabilité portée par $D$. En outre, si $\mu$ est portée par $D$, alors on peut choisir la suite $(X_n)_{n \geqslant 1}$ à valeurs dans $D$ : il suffit de considérer l'inclusion
$$
(X_n)_{n \geqslant 1} : \big( D^{\mathbb{N}^*}, (\mathcal{F}|_D)^{\otimes \mathbb{N}^*}, (\mu|_D)^{\otimes \mathbb{N}^*} \big) \hookrightarrow (X^{\mathbb{N}^*}, \mathcal{F}^{\otimes \mathbb{N}^*})
$$
Ainsi, si $D$ est convexe, $\overline{X}_n$ est également à valeurs dans $D$. À chaque fois qu'on considérera une probabilité $\mu$ portée par une partie convexe $D$ de $X$, les variables aléatoires $X_n$ associées seront toutes supposées à valeurs dans $D$.

\section{Théorie de Cramér}

\subsection{\'Enoncé des résultats principaux}

Soient $(X, \mathcal{C}_0, \mathcal{F}, \tau)$ un e.v.l.c.m., $\mu$ une mesure de probabilité sur $X$ et $(\mu_n)_{n \geqslant 1}$ la suite de Cramér associée à $\mu$. On rappelle la notation
$$
\mathcal{C}_x = \{ x + C \, ; \, C \in \mathcal{C}_0 \}
$$
Le fait que $\mathcal{C}_x$ soit un système fondamental de voisinages de $x$ montre que l'entropie (qui est bien définie ainsi ; cf. Questions de mesurabilité) a pour expression
$$
\forall x \in X \qquad s(x) = \inf_{C \in \mathcal{C}_x} \liminf_{n \to \infty} \frac{1}{n} \log \mathbb{P}\big( \overline{X}_n \in C \big)
$$
La pression est également bien définie (cf. Questions de mesurabilité) et sa fonction convexe-conjuguée a pour expression
$$
\forall x \in X \qquad p^*(x) := \sup_{\lambda \in X^*}  \Big( \langle \lambda | x \rangle - \log \mathbb{E} \big( e^{\langle \lambda | X_1 \rangle} \big) \Big)
$$
On définit le support de $\mu$ par
$$
\textup{supp}(\mu) = \{ x \in X \, ; \, \forall C \in \mathcal{C}_x \quad \mu(C) > 0 \}
$$
et on note $\textup{cosupp}(\mu)$ l'enveloppe fermée convexe du support de $\mu$, \emph{i.e.} l'intersection de tous les ensembles convexes fermés contenant $\textup{supp}(\mu)$. Si $D$ est une partie de $X$ et $\nu$ une probabilité sur $X$ portée par $D$, on note $\nu|_D$ la probabilité trace de $\nu$ sur $D$. On définit enfin, pour $f : X \rightarrow [-\infty , +\infty]$, le domaine de $f$ par
$$
\textup{dom}(f) = \{ x \in X \, ; \, f(x) \in \mathbb{R} \}
$$

\begin{Th}
Soient $(X, \mathcal{C}_0, \mathcal{F}, \tau)$ un e.v.l.c.m., $D$ une partie convexe de $X$ et $\mu$ une probabilité sur $X$ portée par $D$. Soit $(\mu_n)_{n \geqslant 1}$ la suite de Cramér associée à $\mu$. Alors\\
$\bullet$ la suite $(\mu_n)_{n \geqslant 1}$ vérifie \textup{\textsf{(BS$_{\flat , D}$)}} ; en particulier, $(\mu_n)_{n \geqslant 1}$ et $(\mu_n|_D)_{n \geqslant 1}$ vérifient un PGD faible ; et, si $D$ est relativement compact, $(\mu_n)_{n \geqslant 1}$ et $(\mu_n|_D)_{n \geqslant 1}$ vérifient un PGD ;\\
$\bullet$ le domaine de $s$ et l'enveloppe convexe fermée du support de $\mu$ sont contenus dans l'adhérence de $D$, \emph{i.e.}
$$
\textup{dom}(s) \subset \overline{D} \qquad \textrm{et} \qquad \textup{cosupp}(\mu) \subset \overline{D}
$$
$\bullet$ l'entropie est égale à l'opposée de la convexe-conjuguée de la pression, \emph{i.e.}
$$
s = -p^*
$$
\end{Th}

La preuve des deux premiers points repose sur les arguments généraux de grandes déviations. L'amélioration principale à ce sujet, par rapport aux textes classiques, est la résolution des problèmes de mesurabilité. Quant à l'égalité $s = -p^*$ vraie en toute généralité, elle repose sur la notion de convexe-tension relativement à des topologies bien choisies sur $X$. Si $D$ est une partie de $X$ et $\mu$ une probabilité sur $X$, on dit que $\mu$ est \emph{convexe-tendue sur $D$} s'il existe une suite $(K_m)_{m \geqslant 1}$ de parties mesurables de $X$ telle que, pour tout $m \geqslant 1$, $K_m \cap D$ soit convexe et relativement compact et vérifiant
$$
\lim_{m\to\infty}\mu(K_m) = 1
$$
Si la mesure $\mu$ est elle-même convexe-tendue, on a de plus les résultats suivants :
\begin{Th}
Soient $(X, \mathcal{C}_0, \mathcal{F}, \tau)$ un e.v.l.c.m., $D$ une partie convexe de $X$ et $\mu$ une probabilité sur $X$ portée par $D$. Soit $(\mu_n)_{n \geqslant 1}$ la suite de Cramér associée à $\mu$. Si $\mu$ est convexe-tendue sur $D$, alors\\
$\bullet$ la suite $(\mu_n)_{n \geqslant 1}$ vérifie \textup{\textsf{(BS$_{c, D}$)}} ;\\
$\bullet$ l'adhérence du domaine de $s$ est l'enveloppe convexe fermée du support de $\mu$, \emph{i.e.}
$$
\overline{\textup{dom}(s)} = \textup{cosupp}(\mu)
$$
\end{Th}

La convexe-tension s'avère être un outil pertinent pour passer de la borne supérieure pour les convexes compacts à la borne supérieure pour les convexes. En particulier, les demi-espaces ouverts vérifient la borne supérieure, ce qui entraîne l'égalité $p = (-s)^*$. Puis, l'égalité duale $s = -p^*$ découle du théorème de Hahn-Banach dans l'espace vectoriel localement convexe $(X, \tau)$. La notion de convexe-tension a déjà été introduite par \cite{BaZ79} et reprise par \cite{Cer07}, mais nous avons simplifié notablement son utilisation (cf. partie sur la convexe-tension).




\medskip

\textbf{Remarque :} Voici un exemple où $s = -p^*$ alors que $\mu$ n'est pas convexe-tendue elle-même. Sur $X = \mathbb{R}^\mathbb{R}$ muni de la tribu cylindrique et de la topologie produit, le premier théorème assure l'égalité $s = -p^*$. Toutefois, si $\mu$ désigne la mesure image de la mesure de Lebesgue par l'application $x \in [0, 1] \mapsto 1_{[x, x+1]}$, $\mu$ n'est pas convexe-tendue.


\subsection{Exemples d'applications}

Voici les cas les plus importants de la théorie de Cramér que contient ce nouveau cadre :

$\bullet$ \textbf{Théorème de Cramér dans un espace de Banach séparable :} Soient $(X, \tau)$ un espace de Banach séparable, $B$ sa boule unité ouverte, $\widetilde{\mathcal{C}}_0 = \{ B \}$ (cf. la section introduisant les e.v.l.c.m.) et $(X, \mathcal{C}_0, \mathcal{F}, \tau)$ l'e.v.l.c.m. associé. Dans ce cas, $\mathcal{C}_0$ est l'ensemble des boules ouvertes de $X$ centrées en $0$ et $\mathcal{F}$ la tribu borélienne sur $X$. Alors, toute suite i.i.d. $(X_n)_{n \geqslant 1}$ sur $X$ vérifie le PGD faible et l'égalité $s = - p^*$. En outre, toute mesure de probabilité sur $X$ est convexe-tendue (cf. convexe-tension). Donc, on a la borne supérieure pour tous les convexes.

$\bullet$ \textbf{Théorème de Cramér en topologie faible :} Soient $X$ et $Y$ deux espaces en dualité, $\widetilde{\mathcal{C}}_0$ (cf. la section introduisant les e.v.l.c.m.) l'ensemble des bandes ouvertes de la forme $\{ | \langle y | \cdot \rangle | < 1 \}$ pour $y \in Y$ et $(X, \mathcal{C}_0, \mathcal{F}, \tau)$ l'e.v.l.c.m. associé. Dans ce cas, $\tau = \sigma(X ; Y)$ est la topologie faible sur $X$, et $\mathcal{F}$ est la tribu cylindrique sur $X$. Alors, toute suite de Cramér $(\mu_n)_{n \geqslant 1}$ sur $X$ vérifie un PGD faible et l'égalité $s = -p^*$.

$\bullet$ \textbf{Théorème de Sanov :} Cas particulier du point précédent, si $(E, \mathcal{E})$ est un espace mesurable, soient $Y$ l'espace vectoriel des fonctions mesurables bornées de $(E, \mathcal{E})$ dans $(\mathbb{R}, \mathcal{B}(\mathbb{R}))$, muni de la norme $||\cdot||_\infty$, et $X = Y^*$. Alors, $X$ et $Y$ sont en dualité et l'ensemble des mesures de probabilité sur $(E, \mathcal{E})$, noté $\mathcal{M}_1^+(E)$, est un sous-ensemble convexe de $X$. Toute suite i.i.d. $(M_n)_{n \geqslant 1}$ de mesures de probabilité sur $(E, \mathcal{E})$ vérifie donc un PGD faible dans $X$ d'entropie $s = -p^*$. La suite $(M_n)_{n \geqslant 1}$ vérifie également un PGD faible dans $\mathcal{M}_1^+(E)$ relativement à la $\tau$-topologie, trace de $\tau$ sur $\mathcal{M}_1^+(E)$. De plus, comme $\mathcal{M}_1^+(E)$ est relativement compact, $(M_n)_{n \geqslant 1}$ vérifie un PGD.

\section{Retour sur les hypothèses et compléments}

\subsection{E.v.l.c.m. : définition équivalente}

Soient $X$ un espace vectoriel réel et
$$
\overleftarrow{\mathcal{X}} = \big( N_i, f_i, f_{ij} \big)_{i \leqslant j}
$$
un système projectif d'espaces vectoriels normés séparables, autrement dit une famille telle que

\medskip

\textsf{(PROJ$_1$)} les indices $i$ et $j$ décrivent un ensemble $(J, \leqslant)$ préordonné filtrant à droite ;

\medskip

\textsf{(PROJ$_2$)} pour tout $i \in J$, $N_i$ est un espace vectoriel normé séparable, $f_i$ une application linéaire de $X$ dans $N_i$ et, pour tout $(i, j) \in J^2$ tel que $i \leqslant j$, $f_{ij}$ est une application linéaire continue de $N_j$ dans $N_i$ ;

\medskip

\textsf{(PROJ$_3$)} pour tout $(i, j, k) \in J^3$ tel que $i \leqslant j \leqslant k$, on a $f_{ii} = id_{N_i}$,
$$
f_i = f_{ij} \circ f_j  \qquad \textrm{et} \qquad f_{ik} = f_{ij} \circ f_{jk}
$$

On note, pour tout $i \in J$, $B_i$ la boule unité ouverte de $N_i$. On vérifie alors que
$$
\mathcal{C}_0 = \{ f_i^{-1}(t_i B_i) \, ; \, i \in J, \, t_i \in ]0, +\infty[ \}
$$
satisfait aux trois axiomes \textup{\textsf{(EVLCM$_1$)}}, \textup{\textsf{(EVLCM$_2$)}} et \textup{\textsf{(EVLCM$_3$)}}. De plus, on vérifie que, si $(X, \mathcal{C}_0, \mathcal{F}, \tau)$ est l'e.v.l.c.m. associé à $\mathcal{C}_0$, alors $\tau$ (resp. $\mathcal{F}$) est la topologie (resp. tribu) initiale pour la famille $(N_i, f_i)_{i \in J}$. En effet, la topologie initiale est définie ainsi. Quant à la tribu, pour tout $i \in I$, si $\mathcal{B}_i$ désigne la tribu borélienne de $N_i$, étant donné que $N_i$ est normé et séparable, $\mathcal{B}_i$ est la tribu engendrée par les $t B_i + u$, pour $t > 0$ et $u \in N_i$. Avec cette approche, on pourra dire que $(X, \mathcal{C}_0, \mathcal{F}, \tau)$ est \emph{l'e.v.l.c.m. associé} au système projectif $\overleftarrow{\mathcal{X}}$. Réciproquement, si $(X, \mathcal{C}_0, \mathcal{F}, \tau)$ est un e.v.l.c.m., pour tout $C \in \mathcal{C}_0$, on peut définir $N_C = X/ \{ x \in X \, ; \, M_C(x) = 0 \}$ l'espace séparé associé à $(X, M_C)$ et $f_C$ la surjection canonique : $N_C$ est un espace vectoriel normé et séparable. De plus, si $B \subset C$, $N_C$ s'injecte canoniquement dans $N_B$ : on note $f_{CB}$ l'injection associée. Alors, $\overleftarrow{\mathcal{X}} = (N_C, f_C, f_{CB})_{B \subset C}$ est un système projectif d'espaces vectoriels normés séparables et $(X, \mathcal{C}_0, \mathcal{F}, \tau)$ est l'e.v.l.c.m. associé à $\overleftarrow{\mathcal{X}}$.

\medskip

\textbf{Remarque :} L'analogue d'une famille $\widetilde{\mathcal{C}}_0$ serait une structure initiale d'espace vectoriels normés séparables. A partir d'une telle famille d'espaces, il est facile de construire un système projectif définissant les mêmes topologie et tribu sur $X$.


\medskip

\textbf{Remarque :} On peut également définir la notion de système projectif d'e.v.l.c.m. avec des applications linéaires mesurables. Il s'avère que la limite projective d'un tel système projectif est encore un e.v.l.c.m.


\subsection{Convexe-tension}

Dans cette section, nous donnons quelques propriétés de la convexe-tension. Soient $X$ un espace vectoriel réel, $\mathcal{F}$ une tribu sur $X$, $\tau$ une topologie séparée sur $X$, $\mu$ une probabilité sur $X$ et $D$ une partie convexe de $X$. On dit que $\mu$ est \emph{convexe-tendue sur $D$} s'il existe une suite $(K_m)_{m \geqslant 1}$ de parties mesurables de $X$ telle que, pour tout $m \geqslant 1$, $K_m \cap D$ soit un convexe relativement compact et vérifiant
$$
\lim_{m\to\infty}\mu(K_m) = 1
$$
Si $A$ est une partie de $X$, on dit que $A$ est convexe (resp. relativement compacte) \emph{sur $D$} si $A \cap D$ est convexe (resp. relativement compacte). Ici, pour tout $m \geqslant 1$, $K_m$ est convexe sur $D$ et relativement compact sur $D$. Conséquence immédiate de la convexe-tension : pour tout $C \in \mathcal{F}$ convexe sur $D$ et pour tout $m \geqslant 1$, $K_m \cap C$ est mesurable, convexe sur $D$, relativement compact sur $D$ et inclus dans $C$ et on a
$$
\lim_{m\to\infty}\mu(K_m \cap C) = \mu(C)
$$

Mentionnons une autre notion introduite dans \cite{BaZ79} : on dit que $\mu$ est \emph{convexe-régulière sur $D$} si, pour tout $C \in \mathcal{F}$ convexe ouvert, il existe une suite $(K_m)_{m \geqslant 1}$ de parties mesurables de $X$, convexes sur $D$ et relativement compactes sur $D$ telle que, pour tout $m \geqslant 1$, $\overline{K_m} \subset C$, et vérifiant
$$
\lim_{m\to\infty}\mu(K_m) = \mu(C)
$$

Voici le lien entre les deux notions :
\begin{Pro}\label{convreg}
Soient $(X, \tau)$ un espace vectoriel localement convexe séparé, $\mathcal{F}$ une tribu sur $X$ stable par translation et dilatation de rapport non nul, $\mu$ une mesure de probabilité sur $X$ et $D$ une partie convexe de $X$. Alors $\mu$ est convexe-régulière sur $D$ si et seulement si $\mu$ est convexe-tendue sur $D$.
\end{Pro}

\textbf{Remarque :} Ce résultat permet de clarifier un point de l'annexe de \cite{BaZ79}. Et, par la même occasion, d'alléger l'hypothèse ($\hat{\textbf{C}}$) introduite par \cite{DeS89} en réponse à cette imprécision, et reprise dans les textes suivants, notamment \cite[6.1.2. (b)]{DeZ93}.

\medskip

\textbf{Démonstration :} L'implication directe est immédiate. Pour la réciproque, soit $C \in \mathcal{F}$ convexe ouvert. Comme $\mathcal{F}$ est invariante par translation, on peut supposer que $0 \in C$. Etant donné que $\mathcal{F}$ est stable par dilatation de rapport non nul, la propriété \ref{convj} permet d'écrire :
$$
\overline{C} = \{ x \in X ; M_C(x) \leqslant 1 \} = \bigcap_{\substack{ r \in \mathbb{Q} \\ r > 1 }} \{ x \in X ; M_C(x) < r \} = \bigcap_{\substack{ r \in \mathbb{Q} \\ r > 1 }} rC
$$
Aussi $\smash{\overline{C}}$ est-il mesurable, de même que ses dilatés de rapports non nuls. Définissons, pour tout $m \geqslant 1$,
$$
C_m = \left( 1 - \frac{1}{m+1} \right) \overline{C}
$$
de sorte que $C_m \subset C$ et $\mu(C_m) \to \mu(C)$. Soit $(K_m)_{m \geqslant 1}$ une suite de parties de $X$ mesurables, convexes sur $D$ et relativement compactes sur $D$ telle que $\mu(K_m) \to 1$. Alors, pour tout $m \geqslant 1$, $C_m \cap K_m$ est mesurable, convexe sur $D$, relativement compact sur $D$, d'adhérence incluse dans $C$ et vérifie
$$
\phantom \qquad\qquad\qquad\qquad \mu(C_m \cap K_m) \geqslant \mu(C_m) - \big(1 - \mu(K_m) \big) \to \mu(C) \qquad\qquad\qquad\qquad \qed
$$

\medskip

Soient maintenant $(X, \mathcal{C}_0, \mathcal{F}, \tau)$ un e.v.l.c.m., $D$ une partie convexe de $X$, $\mu$ une probabilité sur $X$ portée par $D$ et $(\mu_n)_{n \geqslant 1}$ la suite de Cramér associée à $\mu$. Si $\mu$ est convexe-tendue sur $D$ et si $(K_m)_{m \geqslant 1}$ est une suite de parties mesurables, convexes sur $D$ et relativement compactes sur $D$ telle que $\mu(K_m) \to 1$, alors, pour tout $n \geqslant 1$, $\mu_n(K_m) \to 1$ : donc les $\mu_n$ sont toutes convexes-tendues sur $D$ (au passage, la suite $(K_m)_{m \geqslant 1}$ est commune aux $\mu_n$). En effet, si $n \geqslant 1$, on a :
$$
\mu_n(K_m) = \mathbb{P}(\overline{X}_n \in K_m) \geqslant \mathbb{P}(X_1 \in K_m \, ; \, \ldots \, ; \, X_n \in K_m) = \mu(K_m)^n \to 1
$$
Cette remarque permet de simplifier les hypothèses de \cite[chapitre 9]{Cer07}.
On dit qu'un e.v.l.c. $(X, \tau)$ est un espace de Fréchet s'il est complet et si $0$ admet un système fondamentale de voisinages dénombrable (cette seconde condition est équivalente à la métrisabilité).

\begin{Pro}
Si $(X, \tau)$ est un espace de Fréchet séparable et $\mathcal{F}$ sa tribu borélienne, alors toute probabilité sur $X$ est convexe-tendue sur $X$.
\end{Pro}
Les deux arguments de la preuve sont les suivants : toute probabilité sur $X$ est tendue (cf. \cite[théorème 1.3.]{Bil68}) et l'enveloppe convexe fermée d'un compact de $X$ est compacte, d'après le théorème de Krein (cf. \cite[IV.37]{BouEVT}). En particulier, le résultat est vrai sur un espace de Banach séparable muni de sa tribu borélienne.

\section{Questions de mesurabilité}

La séparabilité (ou plutôt la seconde dénombrabilité) est le nerf de la mesurabilité. Le but de cette partie est de montrer que les fonctions $s$ et $p$ sont bien définies, sur $X$ et $X^*$ respectivement. Si $(X, \mathcal{C}_0, \mathcal{F}, \tau)$ est un e.v.l.c.m. et $(X_n)_{n \geqslant 1}$ une suite de variables aléatoires à valeurs dans $X$, pour pouvoir définir l'entropie, on souhaite que les fonctions $\overline{X}_n$ soient des variables aléatoires. Etant donné que $\mathcal{F}$ est stable par dilatation de rapport non nul, la mesurabilité requise est conséquence de la mesurabilité de l'addition vectorielle. Pour ce qui est de la pression, il s'agit de vérifier que toute forme linéaire continue est mesurables. Nous proposons deux démonstrations de ces deux résultats, correspondant aux deux définitions équivalentes des e.v.l.c.m.

\subsection{Version convexes internes}

\begin{Pro}
Soit $(X, \mathcal{C}_0, \mathcal{F}, \tau)$ un e.v.l.c.m. L'addition vectorielle
$$
(x, y) \in \big( X^2, \mathcal{F}^{\otimes 2} \big) \mapsto x+y \in (X, \mathcal{F})
$$
est mesurable.
\end{Pro}
\textbf{Démonstration :} Il suffit de vérifier que, pour tout $C \in \mathcal{C}_0$,
$$
\{ (x, y) \in X^2 ; x + y \in C \} \in \mathcal{F}^{\otimes 2}
$$
Soit $C \in \mathcal{C}_0$. Soit $Q$ une partie de $X$ dénombrable et dense pour la topologie localement convexe $\tau_C$ engendrée par $C$. Montrons que
$$
\{ (x, y) \in X^2 ; x + y \in C \} = \bigcup_{\substack{ u \in Q \\ r \in \mathbb{Q} \cap ]0, 2[ } } \left\{ (x, y) \in X^2 ; \begin{array}{l} M_C(x - u) < r \\ M_C(y + u) < 1 - r \end{array} \right\}
$$
Pour l'inclusion $\subset$, soit $(x, y) \in X^2$ tel que $x+y \in C$. Notons $a = M_C(x+y) < 1$. Soit $r \in \mathbb{Q} \cap ]0, (1-a)/2[$. L'ouvert non vide $x+rC$ contient alors un point $u$ de $Q$. Et on a $M_C(x-u) < r$ (car $C$ est symétrique) et
$$
M_C(y+u) \leqslant M_C(y+x) + M_C(u-x) < a + r < 1 - r
$$
Pour l'inclusion $\supset$, on note que, pour tout $(x, y, u) \in X^3$,
$$
M_C(x+y) \leqslant M_C(x-u) + M_C(y+u)
$$
d'où le résultat en choisissant convenablement $u \in Q$.\qed

\begin{Le}
Soit $X$ un espace vectoriel. Soient $B$ et $C$ deux convexes internes de $X$. Si $B \subset C$ et si la topologie localement convexe $\tau_B$ engendrée par $B$ est séparable, alors $C$ s'écrit comme réunion dénombrable de translatés et de dilatés de rapports non nuls de $B$.
\end{Le}
\textbf{Démonstration :} Soit $Q$ une partie de $X$ dénombrable et dense relativement à $\tau_B$. Montrons que
$$
C \subset \bigcup_{\substack{ u \in Q , \, r \in \mathbb{Q}_+ \\ u + rB \subset C} } u + rB
$$
l'autre inclusion étant immédiate. Soit $x \in C$. Comme $C$ est interne, $M_C(x) < 1$ et il existe $r \in \mathbb{Q}_+$ tel que $0 < r < 1 - M_C(x)$ ; ainsi $x + rB \subset C$. L'ouvert non vide $x + \frac{1}{2}( rB \cap (-rB) )$ relativement à $\tau_B$ contient un point $u$ de $Q$. On a alors
$$
u \in x - \frac{r}{2}B \qquad\textrm{et}\qquad u + \frac{r}{2} B \subset x + rB \subset C
$$
d'où
$$
\phantom \qquad\qquad\qquad\qquad\qquad\qquad\qquad\quad x \in u + \frac{r}{2} B \subset C \qquad\qquad\qquad\qquad\qquad\qquad\qquad\quad \square
$$


\begin{Pro}
Soit $(X, \mathcal{C}_0, \mathcal{F}, \tau)$ un e.v.l.c.m. Toute forme linéaire continue sur $X$ est mesurable.
\end{Pro}
\textbf{Démonstration :} Soit $\lambda \in X^*$. Soit $H_\lambda$ le demi-espace ouvert $\{ x \in X ; \langle \lambda | x \rangle < 1 \}$. Comme $H_\lambda$ est un voisinage de $0$, il existe $C \in \mathcal{C}_0$ tel que $C \subset H_\lambda$. Le lemme précédent montre alors que $H_\lambda$ est réunion dénombrable de translatés et dilatés de rapports non nuls de $C$, donc est mesurable. Comme $\mathcal{F}$ est stable par translation, $\lambda$ est mesurable.\qed

\subsection{Version système projectif}

\begin{Pro}
Soit $(X, \mathcal{C}_0, \mathcal{F}, \tau)$ un e.v.l.c.m. L'addition vectorielle
$$
(x, y) \in \big( X^2, \mathcal{F}^{\otimes 2} \big) \mapsto x+y \in (X, \mathcal{F})
$$
est mesurable.
\end{Pro}
\textbf{Démonstration :} Supposons que $(X, \mathcal{C}_0, \mathcal{F}, \tau)$ soit l'e.v.l.c.m. associé au système projectif $\overleftarrow{\mathcal{X}} = \big( N_i, f_i, f_{ij} \big)_{i \leqslant j}$. Il suffit de montrer que, pour tout $i \in I$, l'application composée $f_i(\cdot + \cdot)$
$$
(x, y) \in \big( X^2, \mathcal{F}^{\otimes 2} \big) \mapsto f_i(x) + f_i(y) \in (N_i, \mathcal{B}_i)
$$
est mesurable. Cette application apparaît comme le chemin du haut du schéma commutatif
$$
\begin{array}{ccc}
(x, y) \in \big( X^2, \mathcal{F}^{\otimes 2} \big) & \longrightarrow & x + y \in (X, \mathcal{F})\\
\downarrow & & \downarrow \\
\big( f_i(x), f_i(y) \big) \in \big( N_i^2, \mathcal{B}_i^{\otimes 2} \big) & \longrightarrow & f_i(x) + f_i(y) \in (N_i, \mathcal{B}_i)
\end{array}
$$
Or, le chemin du bas est mesurable : la flèche verticale est un produit d'applications mesurables, et la flèche horizontale, addition vectorielle dans l'espace à base dénombrable $N_i$, est continue, donc mesurable.\qed

\begin{Pro}\label{dualmes}
Soit $(X, \mathcal{C}_0, \mathcal{F}, \tau)$ un e.v.l.c.m. Toute forme linéaire continue sur $X$ est mesurable.
\end{Pro}
\textbf{Démonstration :} Supposons que $(X, \mathcal{C}_0, \mathcal{F}, \tau)$ soit l'e.v.l.c.m. associé au système projectif $\overleftarrow{\mathcal{X}} = \big( N_i, f_i, f_{ij} \big)_{i \leqslant j}$. Il s'agit simplement de voir que
$$
X^* = \{ \lambda_i \circ f_i \, ; \, i \in J, \, \lambda_i \in N_i^* \}
$$
Soit $\lambda \in X^*$. Comme $\lambda$ est continue, il existe $C \in \mathcal{C}_0$ tel que
$$
C \subset \{ x \in X \, ; \, \langle \lambda | x \rangle < 1 \}
$$
Ecrivons $C$ sous la forme $f_i^{-1}(t_i B_i)$. On voit alors que
$$
\forall x \in X \qquad f_i(x) = 0 \Rightarrow \langle \lambda | x \rangle = 0
$$
En effet, si $\langle \lambda | x \rangle \neq 0$, il existe $t \in \mathbb{R}$ tel que $\langle \lambda | tx \rangle > 1$, donc $f_i(tx) \notin t_i B_i$, d'où $f_i(x) \neq 0$. On en déduit que l'on peut passer au quotient dans
$$
\begin{array}{rcl}
x \in X & \longrightarrow & \langle \lambda | x \rangle \in \mathbb{R} \\
\searrow & & \nearrow \\
 & f_i(x) \in N_i
\end{array}
$$
et l'application quotient est une forme linéaire continue sur $N_i$.\qed

\section{Démonstration des résultats principaux}

Tout au long de cette partie, $(X, \mathcal{C}_0, \mathcal{F}, \tau)$ désigne un e.v.l.c.m. et, pour tout $x \in X$,
$$
\mathcal{C}_x = \{ x + C \, ; \, C \in \mathcal{C}_0 \}
$$
est un système fondamental de voisinages de $x$ ouverts, convexes et mesurables. On se donne aussi $D$ une partie convexe de $X$, $\mu$ une probabilité sur $X$ portée par $D$ et $(\mu_n)_{n \geqslant 1}$ la suite de Cramér associée à $\mu$.

\subsection{Borne inférieure}

Nous avons vu que l'entropie $s$ définie par
$$
\forall x \in X \qquad s(x) = \inf_{C \in \mathcal{C}_x} \liminf_{n \to \infty} \frac{1}{n} \log \mathbb{P}\big( \overline{X}_n \in C \big)
$$
est la plus grande fonction vérifiant la borne inférieure \textsf{(BI)}. Montrons un résultat un peu plus fort. On étend l'addition à $[-\infty , +\infty]$ via
$$
-\infty \,\textup{\d{$+$}}\, (+\infty) = -\infty
$$
De plus, si $f : X \rightarrow [-\infty , +\infty]$, on dit que $f$ est \emph{semi-continue inférieurement} si, pour tout $t \in [-\infty , +\infty]$, l'ensemble
$$
\big\{ x \in X \, ; \, f(x) > t \big\}
$$
est un ouvert de $X$. On note $f_\bullet$ la régularisée semi-continue inférieurement de $f$, autrement dit la plus grande fonction semi-continue inférieurement qui soit inférieure à $f$.
\begin{Th}[Varadhan ouvert]\label{varadouv}
Pour toute fonction $f : X \rightarrow [-\infty , +\infty]$ mesurable,
$$
\liminf_{n \to \infty} \frac{1}{n} \log \mathbb{E} \Big( e^{nf(\overline{X}_n)} \Big) \geqslant \sup_X (f_\bullet \,\textup{\d{$+$}}\, s)
$$
\end{Th}

\textbf{Démonstration :} Soit $f : X \rightarrow [-\infty , +\infty]$ mesurable. Soient $x \in X$, $\delta > 0$ et $M > 0$. Par définition de $f_\bullet$, il existe $C \in \mathcal{C}_x$ tel que, pour tout $y \in C$,
$$
f(y) \geqslant \min \big( f_\bullet(x) - \delta , M \big)
$$
On a alors, pour tout $n \geqslant 1$ :
\begin{align*}
\frac{1}{n} \log \mathbb{E} \Big( e^{nf(\overline{X}_n)} \Big) & \geqslant \frac{1}{n} \log \mathbb{E} \Big( e^{nf(\overline{X}_n)} 1_{\overline{X}_n \in C} \Big)\\
 & \geqslant \min \big( f_\bullet(x) - \delta, M \big) + \frac{1}{n} \log \mathbb{P}\big( \overline{X}_n \in C \big)
\end{align*}
Prenant la limite inférieure en $n$, il vient
$$
\liminf_{n \to \infty}\frac{1}{n} \log \mathbb{E} \Big( e^{nf(\overline{X}_n)} \Big) \geqslant \min \big( f_\bullet(x) - \delta, M \big) + s(x)
$$
On conclut en faisant tendre $\delta$ vers $0$ et $M$ vers $+\infty$, puis en prenant le supremum en $x \in X$.\qed

\medskip

En particulier, si $f : X \rightarrow \mathbb{R}$ est mesurable et $A \in \mathcal{F}$, le résultat précédent appliqué à $f - \delta_A$ (où $\delta_A = +\infty 1_{X \setminus A}$) donne
$$
\liminf_{n \to \infty} \frac{1}{n} \log \mathbb{E} \Big( e^{nf(\overline{X}_n)} 1_{\overline{X}_n \in A} \Big) \geqslant \sup_{\aroi} (f_\bullet \,\textup{\d{+}}\, s)
$$
que l'on peut voir comme une borne inférieure généralisée. D'autre part, on rappelle que la pression est définie (et bien définie ; cf. Questions de mesurabilité) par
$$
\forall \lambda \in X^* \qquad p(\lambda) := \log \mathbb{E} \big( e^{\langle \lambda | X_1 \rangle } \big)
$$

\begin{Cor} \label{spinegiid}
On a l'inégalité
$$
\forall \lambda \in X^* \quad \forall x \in X \qquad p(\lambda) - s(x) \geqslant \langle \lambda | x \rangle
$$
\end{Cor}

\textbf{Démonstration :} Soit $\lambda \in X^*$. Il suffit de remarquer que
$$
p(\lambda) = \liminf_{n \to \infty} \frac{1}{n} \log \mathbb{E} \big( e^{\langle \lambda | \overline{X}_n \rangle} \big)
$$
et d'appliquer le théorème précédent (\ref{varadouv}) à la fonction mesurable et continue $f = \lambda$.\qed

\medskip

\textbf{Remarque :} On obtient aussi l'inégalité équivalente
$$
\forall x \in X \qquad s(x) \leqslant \inf_{\lambda \in X^*} \big( p(\lambda) - \langle \lambda | x \rangle \big) = -p^*(x)
$$
Notons que notre démonstration fait bien apparaître que l'inégalité $p \geqslant (-s)^*$ (et $s \leqslant -p^*$) est une borne inférieure. La démonstration classique de ce point utilise l'inégalité de Tchebychev. Nous la reproduisons ici. Soient $\lambda \in X^*$ et $x \in X$. Définissons, pour $\varepsilon > 0$, le demi-espace ouvert (mesurable)
$$
H = \{ y \in X ; \langle \lambda | y \rangle > \langle \lambda | x \rangle - \varepsilon \}
$$
Si $n \geqslant 1$, l'inégalité de Tchebychev donne
\begin{align*}
\frac{1}{n} \log \mathbb{P} \big( \overline{X}_n \in H \big) & = \frac{1}{n} \log \mathbb{P} \big( n\langle \lambda | \overline{X}_n \rangle - n\langle \lambda | x \rangle + n\varepsilon > 0 \big)\\
 & \leqslant \frac{1}{n} \log \mathbb{E} \big[ \exp \big( n\langle \lambda | \overline{X}_n \rangle - n\langle \lambda | x \rangle + n\varepsilon \big) \big] = p(\lambda) - \langle \lambda | x \rangle + \varepsilon
\end{align*}
car, pour tout $n \geqslant 1$, $\mathbb{E} ( \exp \langle \lambda | n \overline{X}_n \rangle ) = \mathbb{E} ( \exp \langle \lambda | X_1 \rangle )^n$. Etant donné que $H \in \mathcal{C}_x$, en prenant le supremum en $n \geqslant 1$, on obtient $s(x) \leqslant p(\lambda) - \langle \lambda | x \rangle + \varepsilon$. Puis, passant à l'infimum en $\varepsilon > 0$, il vient
$$
p(\lambda) - s(x) \geqslant \langle \lambda | x \rangle
$$

\medskip

\subsection{Lemme sous-additif}

\begin{Le}\label{lsa1iid}
Soit $\big(u(n)\big)_{n\geqslant 1}$ une suite à valeurs dans $[0, +\infty]$. On suppose que

\medskip

\textup{\textsf{(SA)}} $u$ est sous-additive, \emph{i.e.}
$$
\forall m, n \geqslant 1 \qquad u(m+n) \leqslant u(m) + u(n)
$$

Alors,
$$
\liminf_{n\to\infty} \frac{u(n)}{n} = \inf_{n \geqslant 1} \frac{u(n)}{n}
$$
\end{Le}

\textbf{Démonstration :} Par sous-additivité, pour tous $d, m \geqslant 1$,
$$
\frac{u(dm)}{dm} \leqslant \frac{u(m)}{m}
$$
En faisant tendre $d$ vers $\infty$, on obtient
$$
\liminf_{n\to\infty} \frac{u(n)}{n} \leqslant \liminf_{d\to\infty} \frac{u(dm)}{dm} \leqslant \frac{u(m)}{m}
$$
d'où le résultat en passant à l'infimum en $m \geqslant 1$.\qed

\begin{Le}\label{lsa2iid}
Soit $\big(u(n)\big)_{n\geqslant 1}$ une suite à valeurs dans $[0, +\infty]$. On suppose que

\medskip

\textup{\textsf{(SA)}} $u$ est sous-additive ;

\medskip

\textup{\textsf{(C)}} $u$ est contrôlée, \emph{i.e.} il existe $N \geqslant 1$ tel que
$$
\forall n \geqslant N \qquad u(n) < + \infty
$$

Alors, la suite $\big( u(n)/n \big)_{n\geqslant 1}$ converge vers
$$
\inf_{n \geqslant 1} \frac{u(n)}{n}
$$
\end{Le}

\textbf{Démonstration :} Soient $n \geqslant m \geqslant N$. La division euclidienne de $n$ par $m$ s'écrit $n = mq + r$ avec $q \geqslant 1$ et $r \in \{ 0, \ldots , m-1 \}$ ; ainsi, la sous-additivité permet d'écrire
$$
u(n) = u(mq + r) = u\big( m(q-1) + m+r \big) \leqslant (q-1) u(m) + u(m+r)
$$
puis
$$
\frac{u(n)}{n} \leqslant \frac{u(m)}{m} + \frac{1}{n} \max_{0\leqslant i<m} u(m+i)
$$
D'où, comme $u$ est contrôlée, en faisant tendre $n$, puis $m$ vers $\infty$, on obtient
$$
\limsup_{n\to\infty} \frac{u(n)}{n} \leqslant \liminf_{m\to\infty} \frac{u(m)}{m}
$$
autrement dit la suite $\big( u(n)/n \big)_{n\geqslant 1}$ converge. D'après le lemme \ref{lsa1iid}, sa limite est
$$
\phantom \qquad\qquad\qquad\qquad\qquad\qquad\qquad\qquad\quad \inf_{n \geqslant 1} \frac{u(n)}{n} \qquad\qquad\qquad\qquad\qquad\qquad\qquad\qquad\quad \square
$$

On étend l'addition à $[-\infty , +\infty]$ via
$$
-\infty \,\dot{+}\, (+\infty) = +\infty
$$
Soit $f : X \rightarrow [-\infty , +\infty]$. On dit que $f$ est \emph{concave sur $D$} si
$$
\forall (x, y) \in D^2 \quad \forall \alpha \in [0, 1] \qquad f(\alpha x + (1-\alpha)y) \geqslant \alpha f(x) \,\dot{+}\, (1-\alpha) f(y)
$$
Si $A$ est une partie de $X$, on dit que $A$ est convexe (resp. relativement compacte) \emph{sur $D$} si $A \cap D$ est convexe (resp. relativement compacte). La proposition suivante relie la convexité et la sous-additivité :
\begin{Pro}\label{lsames}
Pour toute fonction mesurable $f : X \rightarrow [-\infty , +\infty]$ concave sur $D$,
$$
\limsup_{n \to \infty} \frac{1}{n} \log \mathbb{E}\Big( e^{nf(\overline{X}_n)} \Big) = \sup_{n \geqslant 1} \frac{1}{n} \log \mathbb{E}\Big( e^{nf(\overline{X}_n)} \Big)
$$
\end{Pro}

\textbf{Démonstration :} Soit $f : X \rightarrow [-\infty , +\infty]$ mesurable et concave sur $D$. On vérifie que la suite de terme général
$$
u(n) = - \log \mathbb{E}\Big( e^{nf(\overline{X}_n)} \Big)
$$
est sous-additive. Soient $m, n \geqslant 1$. Comme $f$ est concave sur $D$ et la suite $(X_k)_{1 \leqslant k \leqslant n}$ i.i.d., on a :
$$
\mathbb{E}\Big( e^{(m+n) f(\overline{X}_{m+n})} \Big) \geqslant \mathbb{E}\Big( e^{mf(\overline{X}_m) + nf(\overline{X}_n)} \Big) = \mathbb{E}\Big( e^{mf(\overline{X}_m)} \Big) \mathbb{E}\Big( e^{nf(\overline{X}_n)} \Big)
$$
On conclut à l'aide du lemme \ref{lsa1iid}.\qed

\begin{Pro} \label{convlim}
Pour tout $C \in \mathcal{F}$ convexe sur $D$,
$$
\limsup_{n \to \infty} \frac{1}{n} \log \mathbb{P}\big( \overline{X}_n \in C \big) = \sup_{n \geqslant 1} \frac{1}{n} \log \mathbb{P}\big( \overline{X}_n \in C \big)
$$
Si, de plus, $C$ est le translaté d'un convexe interne (en particulier si $C$ est un convexe ouvert), alors
$$
\lim_{n \to \infty} \frac{1}{n} \log \mathbb{P}\big( \overline{X}_n \in C \big) = \sup_{n \geqslant 1} \frac{1}{n} \log \mathbb{P}\big( \overline{X}_n \in C \big)
$$
\end{Pro}

\textbf{Démonstration :} Soit $C \in \mathcal{F}$ convexe sur $D$. Pour le premier point, il suffit d'appliquer le lemme précédent à la fonction $- \delta_C$ mesurable et concave sur $D$. Supposons maintenant qu'il existe $x \in X$ tel que $C - x$ soit interne. Pour simplifier les notations, on se ramène au cas où $x = 0$. Ou bien $u(n) = +\infty$ pour tout $n \geqslant 1$ auquel cas on a bien la convergence voulue, ou bien il existe $m \geqslant 1$ tel que $u(m) < +\infty$. Dans ce cas, nous allons montrer que la suite $(u(n))_{n \geqslant 1}$ est contrôlée, \emph{i.e.} vérifie \textsf{(C)}. On a :
$$
0 < \mathbb{P}\big( \overline{X}_m \in C \big) = \mathbb{P}\big( M_C(\overline{X}_m) < 1 \big) = \lim_{\varepsilon \to 0} \mathbb{P}\big( M_C(\overline{X}_m) < 1 - \varepsilon \big)
$$
On en déduit qu'il existe $\varepsilon > 0$ tel que
$$
\mathbb{P}\big( M_C(\overline{X}_m) < 1 - \varepsilon \big) > 0
$$
Soit $n \geqslant 1$. Ecrivons la division euclidienne de $n$ par $m$ sous la forme $n = mq + r$ avec $r \in \{ 1, \ldots , m \}$. La convexité et la $\mathbb{R}_+$-linéarité de $M_C$ permettent de montrer que
$$
M_C(\overline{X}_n) \leqslant \frac{m}{n} \left( \sum_{k=0}^{q-1} M_C \Bigg( \frac{1}{m} \sum_{i=mk+1}^{m(k+1)} X_i \Bigg) \right) + \sum_{i=mq+1}^n \frac{M_C(X_i)}{n}
$$
Sachant que $m/n \leqslant 1/q$ et que $(X_k)_{1 \leqslant k \leqslant n}$ est i.i.d., on peut écrire
\begin{align*}
\mathbb{P}\big( M_C(\overline{X}_n) < 1 \big) & \geqslant \mathbb{P}\Bigg( \forall k \in \{ 0, \ldots , q-1 \} \, M_C \Bigg( \frac{1}{m} \sum_{i=mk+1}^{m(k+1)} X_i \Bigg) < 1 - \varepsilon \, ; \\
& \qquad\qquad\qquad\qquad\qquad \forall i \in \{ mq+1, \ldots , n \} \, M_C(X_i) < \frac{n \varepsilon}{r} \Bigg)\\
 & \geqslant \mathbb{P}\Big( M_C(\overline{X}_m) < 1 - \varepsilon \Big)^q \mathbb{P}\Big( M_C(X_1) < \frac{n \varepsilon}{m} \Big)^m
\end{align*}
Comme $M_C$ est finie partout,
$$
1 = \lim_{t \to +\infty} \mathbb{P}\Big( M_C(X_1) < t \Big)
$$
donc il existe $t \in ]0, +\infty[$ tel que $\mathbb{P}( M_C(X_1) < t ) > 0$. Alors,
$$
\forall n \geqslant \frac{mt}{\varepsilon} \qquad \mathbb{P}\big( \overline{X}_n \in C \big) \geqslant \mathbb{P}\big( M_C(\overline{X}_m) < 1 - \varepsilon \big)^q \mathbb{P}( M_C(X_1) < t )^m > 0
$$
donc la suite $(u(n))_{n \geqslant 1}$ est contrôlée, et la conclusion découle du lemme \ref{lsa2iid}.\qed

\medskip

\textbf{Remarque :} En vertu du lemme \ref{lsiid}, le premier point de la proposition précédente est vrai pour $C$ union finie de parties mesurables convexes sur $D$ et le second pour $C$ union finie de convexes internes mesurables.

\subsection{Propriétés de l'entropie}

\begin{Pro}\label{sconciid}
La fonction $s : X \rightarrow [-\infty, 0]$ définie par
$$
s(x) = \inf_{ C \in \mathcal{C}_x } \, \sup_{n \geqslant 1} \frac{1}{n} \log \mathbb{P}\big( \overline{X}_n \in C \big)
$$
est égale à l'entropie de la suite $(\mu_n)_{n \geqslant 1}$. De plus, elle est semi-continue supérieurement et concave.
\end{Pro}

\textbf{Démonstration :} En vertu de la proposition précédente,
$$
s(x) = \inf_{ C \in \mathcal{C}_x } \, \lim_{n \to \infty} \frac{1}{n} \log \mathbb{P}\big( \overline{X}_n \in C \big)
$$
donc $s$ est l'entropie. Montrons la semi-continuité supérieure de $s$. Soient $t \in \mathbb{R}$ et $x \in X$ tels que $s(x) < t$. Par définition de $s(x)$, il existe $C \in \mathcal{C}_0$ tel que
$$
\sup_{n \geqslant 1} \frac{1}{n} \log \mathbb{P}\big( \overline{X}_n \in x + C \big) < t
$$
Alors, pour tout $y \in x + C$, il existe $\varepsilon > 0$ tel que $y + \varepsilon C \subset x + C$, et ainsi
$$
s(y) \leqslant \sup_{n \geqslant 1} \frac{1}{n} \log \mathbb{P}\big( \overline{X}_n \in y + \varepsilon C \big) \leqslant \sup_{n \geqslant 1} \frac{1}{n} \log \mathbb{P}\big( \overline{X}_n \in x + C \big) < t
$$
On en déduit que $s$ est semi-continue supérieurement. Vérifions que $s$ est concave. Soient $x, y \in X$ et $z = (x+y)/2$. Soit $C \in \mathcal{C}_z$. Par continuité de l'addition, il existe $A \in \mathcal{C}_x$ et $B \in \mathcal{C}_y$ tels que $(A+B)/2 \subset C$. Alors
$$
\mathbb{P} \big( \overline{X}_{2n} \in C \big) \geqslant \mathbb{P} \big( \overline{X}_{n} \in A \big) \mathbb{P} \big( \overline{X}_{n} \in B \big)
$$
puis
$$
\lim_{n \to \infty} \frac{1}{2n} \log \mathbb{P}\big( \overline{X}_{2n} \in C \big) \geqslant \lim_{n \to \infty} \frac{1}{2n} \Big( \log \mathbb{P} \big( \overline{X}_{n} \in A \big) + \log \mathbb{P} \big( \overline{X}_{n} \in B \big) \Big) \geqslant \frac{s(x) + s(y)}{2}
$$
On en déduit que $s((x+y)/2) \geqslant (s(x)+s(y))/2$. La concavité de $s$ en découle, sachant que $s$ est semi-continue supérieurement.\qed

\subsection{Borne supérieure faible}

La borne supérieure pour les compacts repose sur l'interversion infimum-supremum qui fait l'objet du lemme suivant.
\begin{Le}\label{lsiid}
Si $\big(u_1(n)\big)_{n\geqslant 1}$, \ldots , $\big(u_r(n)\big)_{n\geqslant 1}$ sont $r$ suites à valeurs dans $[0, +\infty]$, on a l'égalité
$$
\limsup_{n\to\infty} \frac{1}{n} \log \sum_{i=1}^r u_i(n) = \max_{1\leqslant i\leqslant r} \limsup_{n\to\infty} \frac{1}{n} \log u_i(n)
$$
\end{Le}

\textbf{Démonstration :} De l'encadrement (on rappelle que les suites sont positives)
$$
\max_{1\leqslant i\leqslant r} u_i(n) \leqslant \sum_{i=1}^r u_i(n) \leqslant r \max_{1\leqslant i\leqslant r} u_i(n)
$$
on déduit que
$$
\limsup_{n\to\infty} \frac{1}{n} \log \sum_{i=1}^r u_i(n) = \limsup_{n\to\infty} \frac{1}{n} \log \max_{1\leqslant i\leqslant r} u_i(n)
$$
Puis
\begin{eqnarray*}
\limsup_{n\to\infty} \frac{1}{n} \log \max_{1\leqslant i\leqslant r} u_i(n) & = & \lim_{n\to\infty} \sup_{k\geqslant n} \frac{1}{k} \log \max_{1\leqslant i\leqslant r} u_i(k)\\
 & = & \lim_{n\to\infty} \max_{1\leqslant i\leqslant r} \left(\sup_{k\geqslant n} \frac{1}{k} \log u_i(k)\right)\\
 & = & \max_{1\leqslant i\leqslant r} \lim_{n\to\infty} \left(\sup_{k\geqslant n} \frac{1}{k} \log u_i(k)\right)\\
 & = & \max_{1\leqslant i\leqslant r} \limsup_{n\to\infty} \frac{1}{n} \log u_i(n)
\end{eqnarray*}
La troisième égalité vient du fait facile à vérifier que
$$
\max : [-\infty , +\infty]^r \rightarrow [-\infty , +\infty]
$$
est continue.\qed

\medskip

Pour achever la démonstration du PGD faible, on montre un résultat un peu plus fort. On étend l'addition à $[-\infty , +\infty]$ via
$$
-\infty \,\dot{+}\, (+\infty) = +\infty
$$
De plus, si $f : X \rightarrow [-\infty , +\infty]$, on dit que $f$ est \emph{semi-continue supérieurement} si, pour tout $t \in [-\infty , +\infty]$, l'ensemble
$$
\big\{ x \in X \, ; \, f(x) < t \big\}
$$
est un ouvert de $X$. On note $f^\bullet$ la régularisée semi-continue supérieurement de $f$, autrement dit la plus petite fonction semi-continue supérieurement qui soit supérieure à $f$.
\begin{Th}[Varadhan compact]\label{varadcomp}
Pour toute fonction $f : X \rightarrow [-\infty , +\infty]$ mesurable et pour tout $K \in \mathcal{F}$ relativement compact sur $D$,
$$
\limsup_{n \to \infty} \frac{1}{n} \log \mathbb{E} \Big( e^{nf(\overline{X}_n)} 1_{\overline{X}_n \in K} \Big) \leqslant \sup_{\overline{K \cap D}} [f^\bullet \,\dot{+}\, s]
$$
\end{Th}

\textbf{Démonstration :} Soient $K \in \mathcal{F}$ relativement compact sur $D$, $\delta > 0$ et $M < 0$. Par définition de $f^\bullet$, pour tout $x \in X$, il existe $C(x) \in \mathcal{C}_x$ tel que, pour tout $y \in C(x)$,
$$
f(y) \leqslant \max \big( f^\bullet(x) + \delta , M \big)
$$
Par définition de $s(x)$, quitte à réduire un peu $C(x)$, on peut supposer que
$$
\limsup_{n \to \infty} \frac{1}{n} \log \mathbb{P}\big( \overline{X}_n \in C(x) \big) \leqslant \max \big( s(x) + \delta , M \big)
$$
Du recouvrement de $\smash{\overline{K \cap D}}$ par les $C(x)$ avec $x \in \smash{\overline{K \cap D}}$, on peut extraire un sous-recouvrement fini, noté $\{ C(x_i) ; i \in \{ 1, \ldots , r \} \}$. Pour tout $n \geqslant 1$, on a :
\begin{align*}
\frac{1}{n} \log \mathbb{E} \Big( e^{nf(\overline{X}_n)} 1_{\overline{X}_n \in K} \Big) & \leqslant \frac{1}{n} \log \sum_{i=1}^r \mathbb{E} \Big(  e^{n f(\overline{X}_n)} 1_{\overline{X}_n \in C(x_i)} \Big)\\
 & \leqslant \frac{1}{n} \log \sum_{i=1}^r e^{n \max \big( f^\bullet(x_i) + \delta, M \big)} \mathbb{P}\big( \overline{X}_n \in C(x_i) \big)
\end{align*}
Prenant la limite supérieure en $n$ et utilisant le lemme \ref{lsiid}, on obtient :
\begin{align*}
\limsup_{n \to \infty} \frac{1}{n} \log \mathbb{E} \Big( e^{nf(\overline{X}_n)} 1_{\overline{X}_n \in K} \Big)
 & \leqslant \max_{1\leqslant i\leqslant r} \Big( \max \big( f^\bullet(x_i) + \delta , M \big) + \max \big( (s(x_i) + \delta) , M \big) \Big)\\
 & \leqslant \sup_{x \in \smash{\overline{K \cap D}}} \Big( \max \big( f^\bullet(x) + \delta , M \big) + \max \big( (s(x) + \delta) , M \big) \Big)
\end{align*}
ce qui donne le résultat attendu, en faisant tendre $\delta$ vers $0$ et $M$ vers $-\infty$.\qed

\begin{Cor}[Borne supérieure faible]
Pour tout $K \in \mathcal{F}$ relativement compact sur $D$,
$$
\limsup_{n \to \infty} \frac{1}{n} \log \mathbb{P}\big( \overline{X}_n \in K \big) \leqslant \sup_{\overline{K \cap D}} s
$$
En particulier, $(\mu_n)_{n \geqslant 1}$ et $(\mu_n|_D)_{n \geqslant 1}$ vérifient un PGD faible.
\end{Cor}

\textbf{Démonstration :} Il suffit d'appliquer le lemme de Varadhan compact \ref{varadcomp} à $f = 0$ et $K \in \mathcal{F}$ relativement compact sur $D$.\qed

\medskip

Notons que, si $K \in \mathcal{F}$ est convexe relativement compact sur $D$, la proposition \ref{convlim} donne la borne non asymptotique
$$
\sup_{n \geqslant 1} \frac{1}{n} \log \mathbb{P}\big( \overline{X}_n \in K \big) \leqslant \sup_{\overline{K \cap D}} s
$$

\subsection{Pression et entropie}

\begin{Th}\label{psetoile}
On suppose que $\mu$ est convexe-tendue sur $D$. Alors la pression est la fonction convexe-conjuguée de l'opposée de l'entropie, \emph{i.e.}
$$
\forall \lambda \in X^* \qquad p(\lambda) = \sup_{x \in X} \big( \langle \lambda | x \rangle + s(x) \big)
$$
\end{Th}

\textbf{Démonstration :} L'inégalité $p \geqslant (-s)^*$ ne nécessite pas la convexe-tension : elle repose sur l'inégalité de Tchebychev. Il suffit de passer au supremum en $x \in X$ dans \ref{spinegiid} :
$$
\forall \lambda \in X^* \qquad p(\lambda) \geqslant \sup_{x \in X} \big( \langle \lambda | x \rangle + s(x) \big)
$$
C'est l'autre inégalité qui requiert la convexe-tension de la loi $\mu$ de $X_1$. Soit $(K_m)_{m \geqslant 1}$ une suite de parties mesurables de $X$, convexes sur $D$ et relativement compacts sur $D$ telle que $\mu(K_m) \to 1$. L'idée est de conditionner $\overline{X}_n$ dans $K_m$. On remarque tout d'abord que, pour tout $n \geqslant 1$,
\begin{align*}
\log \mathbb{E} \Big( e^{\langle \lambda | X_1 \rangle} 1_{X_1 \in K_m} \Big)& = \frac{1}{n} \log \mathbb{E} \bigg( \exp \Big( n \Big\langle \lambda \Big| \frac{1}{n} (X_1 + \cdots + X_n) \Big\rangle \Big) 1_{X_1 \in K_m} \cdots 1_{X_n \in K_m} \bigg)\\
 & \leqslant \frac{1}{n} \log \mathbb{E} \Big( e^{n \langle \lambda | \overline{X}_n \rangle} 1_{\overline{X}_n \in K_m} \Big)\\
\end{align*}
Or, le lemme de Varadhan compact \ref{varadcomp} donne
$$
\limsup_{n \to \infty} \frac{1}{n} \log \mathbb{E} \Big( e^{n \langle \lambda | \overline{X}_n \rangle} 1_{\overline{X}_n \in K_m} \Big) \leqslant \sup_{x \in \overline{K_m}} \big( \langle \lambda | x \rangle + s(x) \big)
$$
D'où
$$
\log \mathbb{E} \Big( e^{\langle \lambda | X_1 \rangle} 1_{X_1 \in K_m} \Big) \leqslant \sup_{x \in X} \big( \langle \lambda | x \rangle + s(x) \big)
$$
Enfin, en faisant tendre $m$ vers $\infty$, le lemme \ref{ui} ci-après donne
$$
\phantom \qquad\qquad\qquad\qquad\qquad\qquad\quad p(\lambda) \leqslant \sup_{x \in X} \big( \langle \lambda | x \rangle + s(x) \big) \qquad\qquad\qquad\qquad\qquad\qquad\quad \square
$$

\begin{Le}\label{ui}
Soit $g : X \rightarrow [0, +\infty]$ et $(A_m)$ une suite de parties mesurables de $X$ telle que $\mu(A_m) \to 1$. Alors
$$
\mathbb{E}\big( g(X_1) 1_{X_1 \in A_m} \big) \to \mathbb{E}\big( g(X_1) \big)
$$
\end{Le}

\textbf{Démonstration :} Soit $M > 0$. On a
$$
\mathbb{E}\big( g(X_1) 1_{X_1 \in A_m} \big) \geqslant \mathbb{E}\big( (g(X_1) \wedge M) 1_{X_1 \in A_m} \big) \geqslant \mathbb{E}\big( g(X_1) \wedge M \big) - M \big( 1 - \mu(A_m) \big)
$$
Faisant tendre $m$ vers $\infty$, puis $M$ vers $+\infty$, on obtient le lemme.\qed

\subsection{Transformation de Fenchel-Legendre}

Avant de pouvoir passer à l'égalité duale, il nous faut établir un résultat classique d'analyse convexe. De bonnes références à ce sujet sont \cite{Mor67} et \cite{Roc70}. On se donne ici simplement $(X, \tau)$ un espace vectoriel localement convexe. Notons $X^*$ son dual topologique. Définissons, pour $f : X \rightarrow [-\infty, +\infty]$, sa \emph{transformée de Fenchel-Legendre} par
$$
\forall \lambda \in X^* \qquad f^*(\lambda) := \sup_{x \in X}  \big( \langle \lambda | x \rangle - f(\lambda) \big)
$$
La transformée de Fenchel-Legendre de $g : X^* \rightarrow [-\infty, +\infty]$ se définit de façon analogue.
Etendons l'addition à $[-\infty , +\infty ]$ via
$$
\forall a \in \mathbb{R} \qquad \pm\infty \,\textup{\d{$+$}}\, a = \pm\infty \qquad \textrm{et} \qquad -\infty \,\textup{\d{$+$}}\, (+\infty) = -\infty
$$
ainsi que la multiplication par un réel via
$$
\begin{array}{ll} \forall \alpha > 0 & \alpha \cdot (\pm\infty) = \pm\infty \\
\forall \alpha < 0 & \alpha \cdot (\pm\infty) = \mp\infty
\end{array} \qquad \textrm{et} \qquad 0 \cdot (\pm\infty) = 0
$$
On dit que $f : X \rightarrow [-\infty, +\infty]$ est \emph{convexe} si
$$
\forall (x, y) \in X^2 \quad \forall \alpha \in [0, 1] \qquad f(\alpha x + (1-\alpha)y) \leqslant \alpha f(x) \,\textup{\d{$+$}}\, (1-\alpha) f(y)
$$
On notera que la seule fonction convexe prenant la valeur $-\infty$ est la fonction constante de valeur $-\infty$. Les fonctions convexes autres que les deux constantes $\pm \infty$ sont habituellement dénommées \emph{fonctions convexes propres}. On dit que $f : X \rightarrow [-\infty, +\infty]$ est \emph{concave} si $-f$ est convexe. On dit que $q : X \rightarrow [-\infty , +\infty ]$ est \emph{affine} si $q$ est convexe et concave. Les fonctions affines à valeurs dans $[-\infty, +\infty]$ sont alors les fonctions affines habituelles à valeurs dans $]-\infty, +\infty[$ et les deux fonctions constantes $\pm \infty$. Le seul résultat qui nous intéresse ici est le suivant :
\begin{Pro}\label{fliid}
Soit $(X, \tau)$ un espace vectoriel localement convexe et $f : X \rightarrow [-\infty , +\infty ]$. Alors
$$
f^{**} = f
$$
si et seulement si $f$ est convexe et semi-continue inférieurement relativement à la topologie faible $\sigma(X, X^*)$.
\end{Pro}

\textbf{Remarque :} Plus précisément, la transformation de Fenchel-Legendre réalise une bijection des fonctions convexes $\sigma(X, X^*)$-s.c.i. sur $X$ sur les fonctions convexes $\sigma(X^*, X)$-s.c.i. sur $X^*$, de sorte que, si $f$ est convexe et $\sigma(X, X^*)$-s.c.i., sa transformée de Fenchel-Legendre $f^*$ prend le nom de fonction \emph{convexe-conjuguée} de $f$.

\medskip

\textbf{Démonstration :} L'implication directe découle du fait que, pour toute fonction $g : X^* \rightarrow [-\infty , +\infty]$, la fonction $g^*$ est convexe et $\sigma(X, X^*)$-s.c.i. Montrons la réciproque. On remarque que, si $\lambda \in X^*$, $\langle \lambda | \cdot \rangle - f^*(\lambda)$ est la plus grande fonction affine (y compris les deux fonctions affines $\pm \infty$) dirigée par $\lambda$ et inférieure à $f$. Il s'agit donc de voir que $f$ est la borne supérieure de l'ensemble des fonctions affines continue plus petites que $f$ (y compris les deux fonctions affines $\pm \infty$). Si $f = \pm \infty$, le résultat est immédiat. Sinon, $f$ est une fonction convexe propre. Définissons l'épigraphe de $f$ par $epi(f) = \{ (x, t) \in X \times \mathbb{R} ; f(x) \leqslant t \}$. Le fait que $f$ soit convexe et $\sigma(X, X^*)$-s.c.i. assure que $epi(f)$ est convexe et fermé relativement à $\tau$. Soit $(x, t) \in X \times \mathbb{R} \setminus epi(f)$. Le théorème de Hahn-Banach dans l'espace localement convexe $X \times \mathbb{R}$ donne l'existence d'un hyperplan fermé séparant strictement $(x, t)$ et $epi(f)$. Si cet hyperplan n'est pas vertical, il correspond à une fonction affine plus petite que $f$. Sinon, l'hyperplan est de la forme $H \times \mathbb{R}$ où $H$ est un hyperplan affine fermé de $X$ et $f(x) = +\infty$. Soit alors $p$ une fonction affine continue sur $X$ telle que $p(x) > 0$ et $p(y) = 0$ pour tout $y \in H$. Comme $f$ est une fonction convexe propre, il existe une fonction affine continue et finie $q$ inférieure à $f$ : en effet, il existe $x \in X$ tel que $f(x) \in ]-\infty , +\infty[$ et un hyperplan fermé $H$ séparant $(x, f(x) - 1)$ de $epi(f)$ ; cet hyperplan $H$ n'est pas vertical et correspond à une fonction affine continue et finie $q$ inférieure à $f$. Ainsi, pour tout $\alpha > 0$, $q + \alpha p$ est toujours une fonction affine continue inférieure à $f$. Il suffit alors de choisir $\alpha$ tel que $q(x)+\alpha p(x) > t$.\qed

\medskip

On en déduit enfin :
\begin{Th}
Soit $(X, \mathcal{C}_0, \mathcal{F}, \tau)$ un e.v.l.c.m., $\mu$ une probabilité sur $\mathcal{F}$ portée par un convexe $D$, et $(\mu_n)_{n \geqslant 1}$ la suite de Cramér associée. On suppose que $\mu$ est convexe-tendue sur $D$. Alors l'entropie est l'opposée de la fonction convexe-conjuguée de la pression, \emph{i.e.}
$$
\forall x \in X \qquad s(x) = \inf_{\lambda \in X^*} \big( p(\lambda) - \langle \lambda | x \rangle \big)
$$
\end{Th}

\textbf{Démonstration :} Avec les notations de la proposition \ref{fliid}, le théorème \ref{psetoile} s'écrit $p = (-s)^*$. Comme $-s$ est convexe et s.c.i. (cf. la proposition \ref{sconciid}), donc $\sigma(X, X^*)$-s.c.i., la proposition \ref{fliid} donne
$$
\phantom\qquad\qquad\qquad\qquad\qquad\qquad\qquad\, p^* = (-s)^{**} = -s \qquad\qquad\qquad\qquad\qquad\qquad\qquad\,\qed
$$

Passons enfin à l'égalité $s = -p^*$, vraie en toute généralité dans un e.v.l.c.m., sans supposer $\mu$ convexe-tendue.
\begin{Th}
Pour toute probabilité $\mu$ sur un e.v.l.c.m., l'entropie de la suite de Cramér associée est l'opposée de la fonction convexe-conjuguée de la pression, \emph{i.e.}
$$
\forall x \in X \qquad s(x) = \inf_{\lambda \in X^*} \big( p(\lambda) - \langle \lambda | x \rangle \big)
$$
\end{Th}

\textbf{Démonstration :} Elle repose sur le résultat suivant dont on peut trouver la démonstration dans \cite{Pet11c}.

\begin{Th}[Dawson-Gärtner linéaire]
Soient $\smash{\overleftarrow{\mathcal{X}} = (N_i, f_i, f_{ij})_{(i, j) \in J^2}}$ un système projectif d'espaces vectoriels normés séparables et $(X, \mathcal{C}_0, \mathcal{F}, \tau)$ l'e.v.l.c.m. associé. Soient $\mu$ une probabilité sur $X$ et $(\mu_n)_{n \geqslant 1}$ la suite de Cramér associée à $\mu$. Notant $s$ (resp. $s_i$, $p$, $p_i$) l'entropie de $(\mu_n)_{n \geqslant 1}$ (resp. l'entropie de $(\mu_n \circ f_i^{-1})_{n \geqslant 1}$, la pression de $(\mu_n)_{n \geqslant 1}$, la pression de $(\mu_n \circ f_i^{-1})_{n \geqslant 1}$), on a :
$$
s = \inf_{i \in J} (s_i \circ f_i) \quad \textrm{et} \quad - p^* = \inf_{i \in J} (- p_i^* \circ f_i)
$$
En particulier, si, pour tout $i \in J$, $s_i = -p_i^*$, alors $s = -p^*$.
\end{Th}

Il suffit donc de montrer le résultat pour une probabilité $\mu$ sur un espace vectoriel normé séparable $N$ muni de sa tribu borélienne $\mathcal{B}(N)$. Pour ce faire, les théorèmes de Banach-Mazur et de Mazur-Ulam (cf. \cite[chapitre XI]{Ban32}) montrent qu'on peut voir $N$ comme un sous-espace (non nécessairement fermé) de $C = \mathcal{C}([0, 1] ; \mathbb{R})$. Si $\mathcal{B}$ désigne la tribu borélienne de $C$, on a :
$$
\mathcal{B}(N) = \mathcal{B}|_N := \{ A \cap N \, ; \, A \in \mathcal{B} \}
$$
En effet, l'inclusion $N \hookrightarrow C$ est continue, donc mesurable, donc $\mathcal{B}(N) \supset \mathcal{B}|_N$. L'autre inclusion est vient du fait que $\mathcal{B}(N)$ est la tribu engendrée par les ouvert de $N$, autrement dit les ensembles de la forme $U \cap N$ avec $U$ ouvert de $C$ et que $\mathcal{B}|_N$ contient ces ensembles. On définit alors une mesure $\tilde{\mu}$ sur $\mathcal{B}$ par
$$
\forall A \in \mathcal{B} \qquad \tilde{\mu}(A) = \mu(A \cap N)
$$
qui est une probabilité (portée par $N$ ; \emph{cf.} Suites de Cramér et probabilités portées par une partie). Or, comme $C$ est un espace de Banach séparable, $\tilde{\mu}$ est convexe-tendue sur $C$. Pour éviter toute ambiguïté, notons $s_\mu$ l'entropie de la suite de Cramér associée à $\mu$ et, de façon analogue, $p_\mu$, $s_{\tilde{\mu}}$ et $p_{\tilde{\mu}}$. Le théorème précédent montre que $s_{\tilde{\mu}} = -p_{\tilde{\mu}}^*$. Reste à appliquer le théorème de Dawson-Gärtner linéaire pour le système projectif réduit à l'inclusion $N \hookrightarrow C$ pour en déduire $s_\mu = -p_\mu^*$.\qed




\subsection{Borne supérieure pour les convexes}

Dans toute cette section, on suppose que $\mu$ est convexe-tendue sur $D$. On commence par énoncer une généralisation du lemme \ref{psetoile}. On étend l'addition à $[-\infty , +\infty]$ via
$$
-\infty \,\dot{+}\, (+\infty) = +\infty
$$
On rappelle que $f : X \rightarrow [-\infty , +\infty]$ est concave sur $D$ si
$$
\forall (x, y) \in D^2 \quad \forall \alpha \in [0, 1] \qquad f(\alpha x + (1-\alpha)y) \geqslant \alpha f(x) \,\dot{+}\, (1-\alpha) f(y)
$$
et qu'on note
$$
\textup{dom}(f) = \{ x \in X \, ; \, f(x) \in \mathbb{R} \}
$$
Si $f$ est concave sur $D$, $\textup{dom}(f)$ est convexe sur $D$ (on rappelle que, si $f$ prend la valeur $+\infty$, alors $f = +\infty$, donc $\textup{dom}(f) = \emptyset$ ; \emph{cf.} remarques précédant la proposition \ref{fliid} en les adaptant aux fonctions concaves).

\begin{Th}[Varadhan convexe]\label{varadconv}
Pour toute fonction $f : X \rightarrow [-\infty , +\infty]$ mesurable et concave sur $D$,
$$
\forall n \geqslant 1 \qquad \frac{1}{n} \log \mathbb{E} \Big( e^{nf(\overline{X}_n)} \Big) \leqslant \sup_{X} [f^\bullet \,\dot{+}\, s]
$$
Si, de plus, $\textup{dom}(f)$ est un convexe ouvert et $f$ est semi-continue supérieurement sur $\textup{dom}(f) \cap \overline{D}$, alors
$$
\forall n \geqslant 1 \qquad \frac{1}{n} \log \mathbb{E} \Big( e^{nf(\overline{X}_n)} \Big) \leqslant \sup_{X} [f \,\dot{+}\, s]
$$
\end{Th}

\textbf{Démonstration :} Soient $f : X \rightarrow [-\infty , +\infty]$ mesurable et concave sur $D$, et $N \geqslant 1$. Soit $(K_m)_{m \geqslant 1}$ une suite de parties mesurables de $X$, convexes sur $D$ et relativement compacts sur $D$ telle que $\mu_N(K_m) \to 1$. Pour tout $m \geqslant 1$, le lemme \ref{lsames} appliqué à $f - \delta_{K_m}$ donne
$$
\frac{1}{N} \log \mathbb{E} \Big( e^{Nf(\overline{X}_N)} 1_{\overline{X}_N \in K_m} \Big) \leqslant \limsup_{n \to \infty} \frac{1}{n} \log \mathbb{E} \Big( e^{n f(\overline{X}_n)} 1_{\overline{X}_n \in K_m} \Big)
$$
Puis, le lemme de Varadhan compact \ref{varadcomp} donne
$$
\limsup_{n \to \infty} \frac{1}{n} \log \mathbb{E} \Big( e^{n f(\overline{X}_n)} 1_{\overline{X}_n \in K_m} \Big) \leqslant \sup_{\overline{K_m \cap D}} \big( f^\bullet \,\dot{+}\, s \big)
$$
D'où
$$
\frac{1}{N} \log \mathbb{E} \Big( e^{Nf(\overline{X}_N)} 1_{\overline{X}_N \in K_m} \Big) \leqslant \sup_{X} \big( f^\bullet \,\dot{+}\, s \big)
$$
En faisant tendre $m$ vers $\infty$, le lemme \ref{ui} donne
$$
\frac{1}{N} \log \mathbb{E} \Big( e^{Nf(\overline{X}_N)} \Big) \leqslant \sup_{X} \big( f^\bullet \,\dot{+}\, s \big)
$$
Si $C = \textup{dom}(f)$ est un convexe ouvert, on affine la borne supérieure en choisissant $K_m$ de sorte que $\overline{K_m} \subset C$ et que $\mu(K_m) \to \mu(C)$ (cf. \ref{convreg}). Dans ce cas, si $f$ est semi-continue supérieurement sur $C \cap \overline{D}$, alors $f^\bullet = f$ sur $\overline{K_m \cap D}$. On conclut comme précédemment avec le lemme \ref{ui}.\qed

\begin{Cor}
Soit $f : X \rightarrow [-\infty , +\infty]$ une fonction mesurable et concave sur $D$. Si $\textup{dom}(f)$ est un convexe ouvert et si $f$ est continue sur $\textup{dom}(f) \cap \overline{D}$, alors
$$
\lim_{n \to \infty} \frac{1}{n} \log \mathbb{E} \Big( e^{nf(\overline{X}_n)} \Big) = \sup_{X} \big( f \,\dot{+}\, s \big)
$$
\end{Cor}


\textbf{Remarque :} Si $f$ est concave et s.c.i. sur $\textup{dom}(f)$, alors $\textup{dom}(f)$ est un convexe ouvert.

\medskip

\textbf{Démonstration :} Si $f = +\infty$, le résultat est immédiat. Sinon, comme $f$ est concave, $f$ ne prend pas la valeur $+\infty$. Il suffit alors de combiner les lemmes de Varadhan ouvert \ref{varadouv} et convexe \ref{varadconv}.\qed

\begin{Cor}[Borne supérieure pour les convexes]
Pour tout $C \in \mathcal{F}$ convexe sur $D$,
$$
\forall n \geqslant 1 \qquad \frac{1}{n} \log \mu_n(C) \leqslant \sup_{\overline{C}} s
$$
Si, de plus, $C$ est un convexe ouvert,
$$
\lim_{n \to \infty} \frac{1}{n} \log \mu_n(C) = \sup_{C} s
$$
\end{Cor}

\textbf{Remarque :} La première inégalité est l'objet de \cite[lemma 2.6.]{BaZ79}, mentionné sans démonstration. La version affinée est appelée \og{}inégalité de Bernstein\fg{} dans \cite[théorème 2]{Bar78}. 

\textbf{Démonstration :} Il suffit d'appliquer le lemme de Varadhan convexe \ref{varadconv} et son corollaire à la fonction mesurable $f = - \delta_C := - \infty 1_{X \setminus C}$ qui est concave sur $D$.\qed

\subsection{Domaine de $s$}

Dans toute cette section, $\mu$ est une probabilité sur $X$ portée par un convexe $D$. On définit le support de $\mu$ par
$$
\textup{supp}(\mu) = \{ x \in X \, ; \, \forall C \in \mathcal{C}_x \quad \mu(C) > 0 \}
$$
Cette définition est cohérente avec la définition habituelle, car $\mathcal{C}_x$ est une base de voisinages de $x$. On note $\textup{cosupp}(\mu)$ l'enveloppe convexe fermée de $\textup{supp}(\mu)$, \emph{i.e.} l'intersection de tous les ensembles convexes fermés contenant $\textup{supp}(\mu)$. On rappelle enfin que, pour $f : X \rightarrow [-\infty , +\infty]$, le domaine de $f$ est défini par
$$
\textup{dom}(f) = \{ x \in X \, ; \, f(x) \in \mathbb{R} \}
$$
Si $f$ est une fonction concave, alors $\textup{dom}(f)$ est convexe (on rappelle que, si $f$ prend la valeur $+\infty$, alors $f = +\infty$, donc $\textup{dom}(f) = \emptyset$ ; \emph{cf.} remarques précédant la proposition \ref{fliid} en les adaptant aux fonctions concaves).

\begin{Pro}
On a les inclusions :
$$
\textup{dom}(s) \subset \overline{D} \qquad \textrm{et} \qquad \textup{cosupp}(\mu) \subset \overline{D}
$$
\end{Pro}

\textbf{Démonstration :} Soit $x \in X \setminus \overline{D}$. Il existe $C \in \mathcal{C}_x$ tel que $C \cap D = \emptyset$. Pour tout $n \geqslant 1$, comme $\mu_n$ est portée par $D$, $\mu_n(C) = 0$. Donc $s(x) = 0$ et $x \notin \textup{supp}(\mu)$. Comme $D$ est convexe, on en déduit aussi que $\textup{cosupp}(\mu) \subset \overline{D}$.\qed

\medskip

\textbf{Remarque :} On peut affiner le résultat précédent ainsi : on note $A(D)$ l'ensemble des points $x \in X$ tels que, pour tout demi-espace fermé mesurable $H$ contenant $x$, $\mu(H) > 0$. Alors,
$$
\textup{dom}(s) \subset A(D) \subset \overline{D}
$$
La preuve est analogue à celle de \cite[p. 88]{Cer07}.

\medskip

On suppose désormais que $\mu$ est convexe-tendue sur $D$.

\begin{Th}
On a l'égalité
$$
\overline{\textup{dom}(s)} = \textup{cosupp}(\mu)
$$
\end{Th}

Ce théorème complète les résultats de \cite{BaZ79} et de \cite{Cer07}. Sa démonstration repose sur le lemme suivant :
\begin{Le}
Soient $(X, \tau)$ un espace vectoriel localement convexe séparé, $\mathcal{F}$ une tribu sur $X$ stable par translation et dilatation de rapport non nul, $D$ une partie convexe de $X$ et $\mu$ une mesure de probabilité sur $X$ portée par $D$. Si $\mu$ est convexe-tendue sur $D$, alors, pour tout $C$ convexe ouvert mesurable de $X$,
$$
C \cap \textup{supp}(\mu) \neq \emptyset \iff \mu(C) > 0
$$
\end{Le}

\textbf{Démonstration :} L'implication directe découle de la définition du support de $\mu$ (la convexe-tension n'est pas requise ici). Pour ce qui est de la réciproque, supposons $C \cap \textup{supp}(\mu) = \emptyset$. Comme $\mu$ est convexe-régulière sur $D$, il existe une suite $(K_m)_{m \geqslant 1}$ de parties mesurables de $X$, convexes sur $D$ et relativement compactes sur $D$ telle que $\overline{K_m} \subset C$ et $\mu(K_m) \to \mu(C)$. Comme $\overline{K_m \cap D} \subset X \setminus \textup{supp}(\mu)$, pour tout $x \in \overline{K_m \cap D}$, il existe $C(x) \in \mathcal{C}_x$ tel que $\mu(C(x)) = 0$. Après extraction d'un sous-recouvrement fini $\{ C(x_i) \, ; \, i \in \{ 1, \ldots , r \} \}$ de $\overline{K_m \cap D}$, on obtient
$$
\mu(K_m) = \mathbb{P}(X_1 \in K_m) \leqslant \mathbb{P} \big( X_1 \in C(x_1) \cup \cdots \cup C(x_r) \big) \leqslant \sum_{i=1}^r \mu\big(C(x_i) \big) = 0
$$
d'où $\mu(C) = 0$.\qed

\medskip

\textbf{Remarque :} Si $C = \textup{cosupp}(\mu)$ est mesurable et d'intérieur non vide, un raisonnement analogue montre que $\mu(C) = 1$. En effet, supposant que $0 \in \cro$ et notant $C_m = (1+1/m) \cro$, on a $\overline{K_m \cap (X \setminus C_m)} \subset X \setminus C$ et $\mu(K_m \cap (X \setminus C_m)) \to 1$.

\medskip

\textbf{Démonstration du théorème :} Etant donné que $\textup{dom}(s)$ est convexe, il suffit de montrer les deux inclusions :
$$
\textup{supp}(\mu) \subset \overline{\textup{dom}(s)}
$$
et
$$
\textup{dom}(s) \subset \textup{cosupp}(\mu)
$$

$\bullet$ Soit $x \in X \setminus \overline{\textup{dom}(s)}$. Il existe $C \in \mathcal{C}_x$ tel que $C \cap \textup{dom}(s) = \emptyset$. La borne supérieure améliorée pour les convexes ouverts mesurables donne
$$
\log \mu(C) \leqslant \sup_C s = -\infty
$$
Donc $\mu(C) = 0$ et le lemme permet d'obtenir $C \cap \textup{supp}(\mu) = \emptyset$. D'où la première inclusion.

$\bullet$ Soient $x \in \textup{dom}(s)$ et $C \in \mathcal{C}_0$. Montrons que $(x+2C) \cap \textup{cosupp}(\mu) \neq \emptyset$. Il existe $n \geqslant 1$ tel que $\mu_n(x+C) > 0$. Notons alors
$$
\widetilde{C} = \left\{ (x_1, \ldots , x_n) \in X^n \, ; \, \frac{1}{n}(x_1 + \cdots + x_n) \in x + C \right\}
$$
L'ensemble $\widetilde{C}$ est un ouvert de $X^n$. Soit $Q$ une partie dénombrable de $X$ dense pour la topologie localement convexe $\tau(C)$ engendrée par $\{ tC \, ; \, t \in \mathbb{R}^* \}$. Alors, $Q^n \cap \widetilde{C}$ est dense dans $\widetilde{C}$ pour la topologie $\tau(C)^n$ et
$$
\widetilde{C} \subset \bigcup_{(u_1, \ldots , u_n) \in Q^n \cap \widetilde{C}} \prod_{i=1}^n (u_i + C)
$$
Etant donné que l'union ci-dessus est dénombrable et que $\mu^{\otimes n}(\widetilde{C}) = \mu_n(x+C) > 0$, il existe $(u_1, \ldots , u_n) \in \widetilde{C}$ tel que, pour tout $i \in \{ 1, \ldots , n \}$, $\mu(u_i + C) > 0$. Le lemme précédent assure l'existence, pour tout $i \in \{ 1, \ldots , n \}$, de $y_i \in (u_i + C) \cap \textup{supp}(\mu)$. D'où
$$
y = \frac{1}{n}(y_1 + \cdots + y_n) \in \frac{1}{n}\big( (u_1 + C) + \cdots + (u_n + C_n) \big) \subset x + 2C
$$
Comme de plus $y \in \textup{cosupp}(\mu)$, on obtient le résultat voulu.\qed

\medskip

\textbf{Remarque :} Notons $\mu_n$ la loi de $\overline{X}_n$ et introduisons l'ensemble :
$$
S = \bigcup_{n \geqslant 1} \textup{supp}(\mu_n)
$$
On montre que
$$
\overline{S} = \overline{\textup{dom}(s)} = \textup{cosupp}(\mu)
$$

$\bullet$ Pour voir que $\overline{S} \subset \overline{\textup{dom}(s)}$, il suffit de raffiner le premier point de la preuve ci-dessus : pour $x \in X \setminus \overline{\textup{dom}(s)}$, il existe $C \in \mathcal{C}_x$ tel que $C \cap \textup{dom}(s) = \emptyset$. Soit $n \geqslant 1$. La borne supérieure améliorée pour les convexes ouverts mesurables donne
$$
\frac{1}{n} \log \mu_n(C) \leqslant \sup_C s = -\infty
$$
Donc $\mu_n(C) = 0$ et le lemme permet d'obtenir $C \cap \textup{supp}(\mu_n) = \emptyset$.

\medskip

$\bullet$ Montrons que $\textup{cosupp}(\mu) \subset \overline{S}$. Soit $x \in \textup{cosupp}(\mu)$ et $C \in \mathcal{C}_0$. En vertu du lemme, pour montrer que $x \in \overline{S}$, il suffit de montrer qu'il existe $n \geqslant 1$ tel que $\mu_n(x+C) > 0$. Soient $(x_1, \ldots , x_r) \in \textup{supp}(\mu)^r$ et $(\alpha_1, \ldots , \alpha_r) \in [0, 1]^r$ tels que
$$
\sum_{i=1}^r \alpha_i = 1 \qquad \textrm{et} \qquad \sum_{i=1}^r \alpha_i x_i = x
$$
Par continuité des applications d'espace vectoriel, il existe $(U_1, \ldots , U_r)$ tel que, pour tout $i \in \{ 1, \ldots , r \}$, $U_i$ soit un voisinage de $\alpha_i$ pour la topologie induite sur $[0, 1]$ et tel que
$$
\sum_{i=1}^r U_i x_i \subset x + C
$$
Soit maintenant, pour tout $i \in \{ 1, \ldots , r \}$, un rationnel $r_i \in U_i$. Prenant un dénominateur commun $n$, on écrit $r_i = k_i/n$ et on définit
$$
\widetilde{x} = \sum_{i=1}^r \frac{k_i}{n} x_i \in x + C
$$
Il existe alors $\varepsilon > 0$ tel que $\widetilde{x} + \varepsilon C \subset x + C$ et on a :
$$
\phantom\qquad\qquad\qquad\qquad \mu_n(x + C) \geqslant \mu_n(\widetilde{x} + \varepsilon C) \geqslant \prod_{i=1}^r \mu(x_i + \varepsilon C) > 0 \qquad\qquad\qquad\qquad\qed
$$

\newpage

\section{Compléments}

\subsection{Régularité de $p$}

La fonction $p$ est convexe (conséquence de l'inégalité de Jensen) et, dans le cas où $\mu$ est convexe-tendue, $p$ est s.c.i. (cf. \cite[Lemme 12.1.]{Cer07}). On notera que ces propriétés ne nous servent pas au cours de la démonstration.

\subsection{Un exemple où $p \neq (-s)^*$}

On se place sur $X = \ell^\infty = L^\infty(\mathbb{N} ; \mathbb{R})$. On munit $X$ de la topologie $\tau$ associée à la norme $\| \cdot \|_\infty$. Munissons $X$ d'une tribu.
On note $\mathcal{P}^*$ l'ensemble des parties de $\mathbb{N}$ dont le complémentaire est infini, et on définit
$$
S = \{ 1_P \in X \, ; \, P \in \mathcal{P}^* \}
$$
Puis
$$
\pi : \alpha \in [0, 1] \longmapsto (\alpha_k)_{k \in \mathbb{N}} \in S
$$
où $\overline{\alpha_0\alpha_1 \cdots}$ est l'écriture binaire de $\alpha$. Cette écriture est unique si l'on suppose la suite non cofinale à $1$, et c'est le cas vu la définition de $S$. L'application $\pi$ est donc une bijection. On vérifie que
$$
\mathcal{F} = \Big\{ A \subset X ; \, \forall x \in X \quad \pi^{-1}\big( (A+x) \cap S \big) \in \mathcal{B}(\mathbb{R}) \Big\}
$$
est une tribu, invariante par translation et dilatation de rapport non nul, et contenant les boules ouvertes. De plus,
$$
\big\{ (x, y) \in X^2 \, ; \, x + y \in B(0, 1) \big\} = \bigcap_{k \in \mathbb{N}} \big\{ (x, y) \in X^2 \, ; \, |x_k + y_k| < 1 \big\} \in \mathcal{F} \otimes \mathcal{F}
$$
Donc, pour toute boule ouverte $B$ de $X$,
$$
\big\{ (x, y) \in X^2 \, ; \, x + y \in B \big\} \in \mathcal{F} \otimes \mathcal{F}
$$
On munit $X$ de la tribu $\mathcal{F}$. On définit maintenant, sur $\mathcal{F}$, une mesure par
$$
\forall A \in \mathcal{F} \qquad \mu(A) = \textup{Leb} \circ \pi^{-1} (A \cap S)
$$
et on note $(\mu_n){n \geqslant 1}$ la suite de Cramér associée à $\mu$. Bien que l'addition vectorielle ne soit pas mesurable \emph{a priori}, les considérations précédentes permettent de définir l'entropie $s$ associée par
$$
\forall x \in X \qquad s(x) = \inf_{r > 0} \liminf_{n \to \infty} \frac{1}{n} \log \mu_n \big( B(x, r) \big)
$$
Montrons que $s = -\infty$. Calculons $s(0)$, les autres valeurs se calculant de la même manière. Pour cela, montrons que, pour tout $n \geqslant 1$,
$$
\mu_n \big( B(0, 1) \big) = 0
$$
Pour $n = 1$, la définition de $\pi$ donne :
$$
\mu \big\{ x \in X \, ; \, |x_0| < 1, \ldots , |x_k| < 1 \big \} = \textup{Leb} \big[ 0, 2^{-(k+1)} \big[ = 2^{-(k+1)}
$$
donc
$$
\mu \big( B(0, 1) \big) = 0
$$
Pour $n = 2$, on remarque que
$$
\big\{ (x, y) \in X^2 \, ; \, |x_0 + y_0| < 2 \big\} \cap (S \times S) = \Bigg\{ (1_P, 1_Q) \, ; \, (P, Q) \in \mathcal{P}^* \textrm{ et } \left\{ \begin{array}{ll} P \subset \mathbb{N} \setminus \{ 0 \} \\ \textrm{ou } Q \subset \mathbb{N} \setminus \{ 0 \} \end{array} \right. \Bigg\}
$$
d'où
$$
\mu_2 \big\{ (x, y) \in X^2 \, ; \, |x_0 + y_0| < 2 \big\} = \frac{3}{4}
$$
puis
$$
\mu_2 \big\{ (x, y) \in X^2 \, ; \, |x_0 + y_0| < 2, \ldots , |x_k + y_k| < 2 \big\} = \left( \frac{3}{4} \right)^{k+1}
$$
et
$$
\mu_2 \big( B(0, 1) \big) = 0
$$
Le même raisonnement montre que, pour tout $n \geqslant 1$,
$$
\mu_n \big( B(0, 1) \big) = 0
$$
D'où $s(0) = -\infty$, puis $s = -\infty$. Ainsi,
$$
(-s)^*(\lambda) = \sup_{x \in X} \big( \langle \lambda | x \rangle + s(x) \big) = -\infty
$$
Or $p(0) = 0$, donc
$$
p \neq (-s)^*
$$
Toutefois, on ne peut pas prendre les transformées de Fenchel-Legendre car $p$ n'est pas définie sur tout $X^*$, \emph{a priori}\footnote{On remarquera que $\textup{card}(X^*) = 2^{\aleph_1}$ et que la tribu engendrée par les boules ouvertes de $X$ est de cardinal $\aleph_1$. On doit pouvoir montrer que la tribu $\mathcal{F}$ est aussi de cardinal $\aleph_1$.}.

\subsection{Autour de $\textup{cosupp}(\mu)$}

Soit $(X, \mathcal{C}_0, \mathcal{F}, \tau)$ un e.v.l.c.m., $\mu$ une mesure de probabilité sur $X$, et $(\mu_n)_{n \geqslant 1}$ la suite de Cramér associée. Notons
$$
S = \bigcup_{n \geqslant 1} \textrm{supp}(\mu_n)
$$
Pour toute partie $A$ de $X$, on définit
$$
\textup{conv}_\mathbb{Q}(A) = \left\{ x \in X \, ; \, \exists r \geqslant 1 \quad \exists (x_1, \ldots , x_r) \in A^r \quad x = \frac{1}{r}(x_1 + \cdots + x_r) \right\}
$$
On note aussi $\textup{conv}(A)$ l'enveloppe convexe de $A$ et $\textup{co}(A) = \smash{ \overline{\textup{conv}(A)} }$ l'enveloppe convexe fermée de $A$. On montre que $\smash{ \overline{\textup{conv}_\mathbb{Q}(A)} } = \textup{co}(A)$. Il est facile de voir que
$$
\textup{conv}_\mathbb{Q} \big( \textup{supp}(\mu) \big) \subset S
$$
Il n'y a pas d'autre inclusion, en général, entre les ensembles $\textup{conv}_\mathbb{Q}(\textup{supp}(\mu))$, $\textup{conv}(\textup{supp}(\mu))$ et $S$, comme le montrent les exemples suivants :

\medskip

$\bullet$ Sur $X = \mathbb{R}$ muni de sa topologie usuelle et de la tribu borélienne, soit $\mu = (\delta_0 + \delta_1)/2$. Alors
$$
S = \mathbb{Q} \cap [0, 1] = \textup{conv}_\mathbb{Q} \big( \textup{supp}(\mu) \big)
$$
de sorte que $\textup{conv}(\textup{supp}(\mu)) = [0, 1] \nsubseteq S$.

\medskip

$\bullet$ Voici un exemple de probabilité $\mu$ pour laquelle $S \nsubseteq \textup{conv}(\textup{supp}(\mu))$. On se place sur $X = \mathbb{R}^2$ muni de sa topologie usuelle $\tau$ et de la tribu borélienne $\mathcal{F}$. On définit
$$
f : t \in \mathbb{R}^* \mapsto \bigg( t, \frac{1}{|t|} \bigg)
$$
et on considère $\mu$ la mesure image par $f$ de la loi $\mathcal{N}(0, 1)$. Alors, $\textup{conv}(\textup{supp}(\mu))$ est le demi-espace supérieur ouvert $\{ (x, y) \in \mathbb{R}^2 \, ; \, y > 0 \}$ et $\textup{supp}(\mu_2)$ contient l'axe des abscisses $\{ (x, y) \in \mathbb{R}^2 \, ; \, y = 0 \}$.

\medskip

$\bullet$ On notera que, dans $\mathbb{R}$, l'enveloppe convexe d'un fermé est fermée (ce qui n'est plus vrai dès la dimension $2$ --- cf. l'exemple précédent), de sorte que, dans $\mathbb{R}$,
$$
\textup{conv}_\mathbb{Q} \big( \textup{supp}(\mu) \big) \subset S \subset \textup{conv}\big( \textup{supp}(\mu) \big)
$$
Voici un exemple de probabilité $\mu$ sur $\mathbb{R}$ telle que $S \neq \textup{conv}_\mathbb{Q}(\textup{supp}(\mu))$. Il suffit de prendre une mesure $\mu$ de support $\mathbb{Z}$ (par exemple, $\mu(\{ 0 \}) = 1/2$ et, pour tout $n \geqslant 1$, $\mu(\{  n \}) = \mu(\{  - n \}) = 1/2^{-(n+2)}$). Alors $S = \mathbb{R}$ et $\textup{conv}_\mathbb{Q}(\textup{supp}(\mu)) = \mathbb{Q}$.

\subsection{Pour approfondir la notion de convexe interne}

Les outils utilisé s'inspirent principalement de \cite{BouEVT}. On trouvera aussi une autre présentation dans \cite{RoR64}. Soient $C$ un convexe de $X$ et $x \in C$. On dit que $x$ est un \emph{point interne} de $C$ si la jauge $M_{C-x}$ est finie partout, autrement dit si
$$
X = \bigcup_{t \in [0, +\infty[} x + t(C-x)
$$
Un convexe $C$ de $X$ contenant $0$ est \emph{interne} si et seulement si tout point de $C$ est interne. On vérifie que, si $0$ est un point interne de $C$, alors
$$
\{ x \in X \, ; \, M_C(x) < 1 \}
$$
est l'ensemble des points internes de $C$. On dit qu'un convexe $C$ de $X$ est \emph{absorbant} si $0$ est un point interne de $C$. Plus généralement que celle de convexe interne, la notion de convexe absorbant est fondamentale dès que l'on parle de topologie localement convexe : en effet, si $\mathcal{C}_0$ est une famille de convexes absorbants de $X$, stable par intersection finie et par dilatation de rapport non nul, alors $\mathcal{C}_0$ est un système fondamental de voisinages de $0$ pour une topologie localement convexe $\tau$ sur $X$ (cf. \cite[II.25]{BouEVT}). Si de plus les convexes de $\mathcal{C}_0$ sont internes, alors $\mathcal{C}_0$ est un système fondamental de voisinages de $0$ à la fois convexes et ouverts ; et c'est cette propriété supplémentaire qui nous sert dans ce chapitre (car nous avons la limite pour les convexes ouverts mesurables). Notons encore que l'ensemble des convexes internes est l'ensemble des convexes qui sont des voisinages ouverts de $0$ pour une certaine topologie localement convexe sur $X$. Ou plutôt : un convexe $C$ est interne si et seulement si il existe une topologie localement convexe pour laquelle $C$ est un voisinage ouvert de $0$.

\medskip

\label{sepabs}\textbf{Remarque :} Soit $C$ un convexe contenant $0$. Si la topologie engendrée par les translatés et dilatés de rapports non nuls de $C$ est séparable, alors $C$ est absorbant.

\medskip

Terminons cette section par les liens algébriques entre convexes internes et convexes absorbants. Si $C$ est un convexe de $X$ contenant $0$, on peut montrer que l'ensemble des points internes de $C$ est ou bien $\emptyset$, ou bien un convexe interne. Un convexe est donc absorbant si et seulement si il contient un convexe interne. Dans l'autre sens, si $C$ est un convexe de $X$, $C$ est interne si et seulement si $C$ est absorbant et, pour tout $x \in C$,
$$
C = \bigcup_{t \in ]0, 1[} x + t(C-x)
$$
On a vu qu'un ouvert convexe contenant $0$ est interne. Plus généralement, pour tout convexe $C$ de $X$, l'intérieur pour $\tau$ de $C$ est ou bien $\emptyset$, ou bien l'ensemble des points internes de $C$.


\subsection{Limite pour les convexes}

Soient $(X, \mathcal{C}_0, \mathcal{F}, \tau)$ un e.v.l.c.m., $\mu$ une probabilité sur $\mathcal{F}$ et $(\mu_n)_{n \geqslant 1}$ la suite de Cramér associée à $\mu$. On s'intéresse à la question : pour quels ensembles $A$ la suite de terme général
$$
\frac{1}{n} \log \mu_n(A)
$$
est-elle convergente ? La proposition \ref{convlim} montre que, de façon générale, la suite converge dès que $A$ est le translaté d'un convexe interne. Sans cette hypothèse, la suite peut ne pas converger, même pour un convexe : il suffit de considérer, sur $\mathbb{R}$, $\mu = (\delta_{-1} + \delta_{1})/2$ et $A = \{ 0 \}$, de sorte que
$$
\liminf_{n \to \infty} \frac{1}{n} \log \mu_n(A) = -\infty
$$
(valeurs des termes impairs de la suite) et
$$
\limsup_{n \to \infty} \frac{1}{n} \log \mu_n(A) = 0
$$
De manière générale (même en dimension infinie), si $C$ est un convexe mesurable, la démonstration du lemme sous-additif montre qu'il existe $k_C \in \mathbb{N}^* \cup \{ \infty \}$ tel que
$$
\limsup_{n \to \infty} \frac{1}{n} \log \mu_n(C) = \lim_{\substack{n \to \infty \\ n \in k_C \mathbb{N} }} \frac{1}{n} \log \mu_n(C) \in ]-\infty , 0]
$$
et
$$
\forall n \notin k_C \mathbb{N} \qquad \frac{1}{n} \log \mu_n(C) = -\infty
$$
Si $k_C = 1$, la suite converge dans $]-\infty , 0]$ et, si $k_C = \infty$, la suite est constamment égale à $-\infty$. Notons $\textup{cosupp}(\mu)$ l'enveloppe convexe fermée du support de $\mu$, défini par
$$
\textup{supp}(\mu) = \{ x \in X \, ; \, \forall C \in \mathcal{C}_x \quad \mu(C) > 0 \}
$$
On peut adapter la démonstration de \ref{convlim} pour montrer que la suite de terme général
$$
\frac{1}{n} \log \mathbb{P}\big( \overline{X}_n \in C \big)
$$
converge vers son supremum plus généralement :
\begin{itemize}
\item si $C$ est un convexe mesurable dont l'ensemble des points internes rencontre $\textup{cosupp}(\mu)$ ;
\item ou si $C \in \mathcal{F}$ ne rencontre pas $\textup{cosupp}(\mu)$ (auquel cas le supremum vaut $-\infty$).
\end{itemize}


On en déduit le résultat suivant :
\begin{Pro}
Soient $(X, \mathcal{C}_0, \mathcal{F}, \tau)$ un e.v.l.c.m., $\mu$ une probabilité sur $\mathcal{F}$ et $(\mu_n)_{n \geqslant 1}$ la suite de Cramér associée à $\mu$. Si $\mu$ ne charge pas les hyperplans affines de $X$, alors, pour tout convexe $C \in \mathcal{F}$ admettant un point interne, la suite de terme général
$$
\frac{1}{n} \log \mu_n(C)
$$
converge dans $[-\infty , 0]$ vers son supremum.
\end{Pro}

\textbf{Remarque :} Le théorème de Fubini permet de voir que, si $\mu$ ne charge par les hyperplans, il en est de même de $\mu_n$, pour tout $n \geqslant 1$.

\medskip

\textbf{Remarque :} Si le sous-espace affine $\textup{aff}(C)$ engendré par $C$ est un sous-espace affine strict de $X$, comme, pour tout $n \geqslant 1$, $\mu_n$ ne charge pas les hyperplans, la limite existe pour $C$, et vaut $-\infty$. Le cas intéressant est donc seulement celui où $\textup{aff}(C) = X$. Or, dans ce cas, il n'est pas automatique que $C$ admette un point interne. C'est le cas si $X$ est de dimension finie, mais, par exemple, dans $X = \ell^2 = L^2(\mathbb{N} ; \mathbb{R})$,
$$
K = \bigg\{ x \in \ell^2 \, ; \forall n \quad u_n \geqslant 0 \textrm{ et } \sum_{n \in \mathbb{N}} n u_n \leqslant 1 \bigg\}
$$
est un compact convexe (comme enveloppe convexe fermée du compact $\{ 0 \} \cup \{ e_k/k \, ; \, k \geqslant 1 \}$, où $e_k(n) = \delta_{n,k}$, dans l'espace complet séparable $\ell^2$) dont le sous-espace affine engendré est $X$ et qui n'admet pas de point interne. En effet, si $u \in K$ était un point interne, soit $(n_k)_{k \geqslant 1}$ une suite d'entiers strictement croissante telle que, pour tout $k \geqslant 1$,
$$
\sum_{n \geqslant n_k} nu_n \leqslant 3^{-k}
$$
En particulier,
$$
u_{n_k} \leqslant \frac{1}{n_k 3^k}
$$
Définissons, pour tout $n \in \mathbb{N}$,
$$
w_n = \left\{ \begin{array}{ll} -\frac{1}{n_k 2^k} & \textrm{si $n = n_k$}\\ 0 & \textrm{sinon} \end{array} \right.
$$
Alors $w \in \ell^2$ mais il n'existe pas $v \in K$ et $t \geqslant 0$ tels que $w = u + t(v - u)$. En effet, il faut certainement $t \geqslant 1$ et on a :
$$
u_{n_k} + t(v_{n_k} - u_{n_k}) \geqslant (1 - t) u_{n_k} \geqslant \frac{1-t}{n_k 3^k} > -\frac{1}{n_k2^k}
$$
pour $k$ assez grand.

\medskip

\textbf{Démonstration :} Soit $C \in \mathcal{F}$ convexe. Pour simplifier les notations, supposons que $0$ est un point interne de $C$. Dans ce cas, l'ensemble des points internes de $C$ est $\{ M_C < 1 \}$. En vertu de la remarque précédente, il reste donc à traiter le cas où
$$
C \cap \textup{cosupp}(\mu) \neq \emptyset \qquad \textrm{et} \qquad \{ M_C < 1 \} \cap \textup{cosupp}(\mu) = \emptyset
$$
Dans ce cas, $C \cap \textup{cosupp}(\mu)$ est un convexe inclus dans $\{ M_C = 1 \}$ : il est donc inclus dans un hyperplan affine $H$. En effet, $0 \notin \textup{aff}(\{ M_C = 1 \})$, sinon il existerait $x, y \in \{ M_C = 1 \}$ tels que
$$
0 = x + t(y-x)
$$
et on en déduirait que
$$
M_C(y) \leqslant 1 - \frac{1}{t} < 1
$$
Comme, pour tout $n \geqslant 1$, $\mu_n(H) = 0$ (notons que les hyperplans sont mesurables, car $X$ est un e.v.l.c.m.), on en déduit que $s(C) = -\infty$ et on a la limite souhaitée.\qed

\medskip

En dimension finie, on peut affiner le résultat. On se donne désormais un entier positif $d$ et on se place dans $X = \mathbb{R}^d$, muni de la topologie standard et de la tribu borélienne. On se donne une probabilité $\mu$ sur $X$ et on note $(\mu_n)_{n \geqslant 1}$ la suite de Cramér associée.

\begin{Pro}
Soient $\mu$ une mesure de probabilité sur un espace de dimension finie et $(\mu_n)_{n \geqslant 1}$ la suite de Cramér associée. Si, pour tout couple d'hyperplans parallèles et distincts $(F, G)$,
$$
\mu(F)\mu(G) = 0
$$
alors, pour tout convexe mesurable $C$, la suite
$$
\frac{1}{n} \log \mu_n(C)
$$
converge dans $[-\infty , 0]$ vers son supremum.
\end{Pro}

\textbf{Remarque :} Même en dimension finie, un convexe n'est pas automatiquement borélien (pour peu que sa frontière ne le soit pas). En revanche, comme sa frontière est un ensemble négligeable pour la mesure de Lebesgue, tout convexe est un lebesguien. Comme nous considérons ici la tribu borélienne, nous précisons que $C$ est mesurable.

\medskip

\textbf{Démonstration :} On montre le résultat plus généralement pour une mesure $\mu$ de masse inférieure à $1$. La démonstration se fait par récurrence sur la dimension $d$ de l'espace affine ambiant. Le résultat est immédiat pour $d = 0$. Supposons le résultat vrai pour tout entier strictement inférieur à $d$ et soit $\mu$ une mesure finie sur $\mathbb{R}^d$ vérifiant la condition de l'énoncé.
Soit $C$ un convexe mesurable de $\mathbb{R}^d$. Remarquons que, pour tout $n \geqslant 1$, comme $\mu_n$ est régulière,
$$
\mu_n(C) = \mu_n\big( C \cap \textup{cosupp}(\mu) \big)
$$
On distingue trois cas.

\medskip

$\bullet$ Si $C \cap \textup{cosupp}(\mu) = \emptyset$, alors
$$
\lim_{n \to \infty} \frac{1}{n} \log \mu_n(C) = -\infty
$$

$\bullet$ Si $\cro \cap \textup{cosupp}(\mu) \neq \emptyset$, le lemme sous-additif pour les convexes internes montre que
$$
\lim_{n \to \infty} \frac{1}{n} \log \mu_n(C) = \sup_{n \to \infty} \frac{1}{n} \log \mu_n(C) \in \mathbb{R}
$$
On notera que, $\mu$ étant de masse inférieure à $1$, la suite $u(n) = -\log \mu_n(C)$ est bien positive.

\medskip

$\bullet$ Reste le cas où $C \cap \textup{cosupp}(\mu) \neq \emptyset$ et $\cro \cap \textup{cosupp}(\mu) = \emptyset$. Dans ce cas,
$$
\overbrace{C \cap \textup{cosupp}(\mu)}^o = \cro \cap \overbrace{\textup{cosupp}(\mu)}^o = \emptyset
$$
Ainsi, $\widetilde{C} = C \cap \textup{cosupp}(\mu)$, convexe non vide d'intérieur vide, engendre un sous-espace affine strict $\widetilde{X}$ de $X$. Notons $\tilde{\mu} = \mu|_{\widetilde{X}}$ et $(\tilde{\mu}_n)_{n \geqslant 1}$ la suite de Cramér associée. La condition sur $\mu$ assure que, pour tout $n \geqslant 1$,
$$
\mu_n \big( \widetilde{X} \big) = \tilde{\mu}_n \big( \widetilde{X} \big)
$$
(le reste de $\mu_n$ \og{}ne voit pas\fg{} $\widetilde{X}$, qui est contenu dans un hyperplan). De plus, $\tilde{\mu}$ vérifie l'hypothèse de l'énoncé. On conclut en appliquant l'hypothèse de récurrence à $\widetilde{X}$, $\tilde{\mu}$ et $\widetilde{C}$.\qed

\medskip

La réciproque n'est pas vraie : dans $\mathbb{R}$, il suffit de considérer
$$
\mu = \sum_{n=1}^\infty 2^{-n} \delta_n
$$
Dans ce cas, pour tout $a \in \mathbb{Q} \cap [1, +\infty[$, $k_{\{ a \}} = 1$ et, pour tout $a \in \mathbb{R} \setminus (\mathbb{Q} \cap [1, +\infty[)$, $k_{\{ a \}} = \infty$. Le cas des convexes $I$ d'intérieur non vide se traite en distinguant les cas $I \subset ]-\infty , 0[$ ($k_I = \infty$) et $I \cap [0, +\infty[ \neq \emptyset$ ($k_I = 1$).

\newpage

\section{Suppléments techniques}

\subsection{De l'utilité des hypothèses}

La structure d'espace vectoriel localement convexe est naturelle en théorie de Cramér. Ce qui fait fonctionner la sous-additivité, c'est la structure de système fondamental de voisinages convexes (non nécessairement absorbants, \emph{a priori}). La topologie d'espace vectoriel localement convexe (donc les ouverts convexes absorbants) est imposée d'abord car on a, \emph{a priori}, la limite uniquement pour $C$ translaté de convexe interne. D'autre part, la remarque de la page \pageref{sepabs} montre que la séparabilité impose aux convexes d'être absorbants.
La question de l'utilité de la séparabilité pour la mesurabilité est une question plus délicate que nous n'avons pas plus abordé que la remarque faite en introduction.

\subsection{Quelques extensions du cadre}

Soient $X$ un espace vectoriel, $\mathcal{F}$ une tribu sur $X$ stable par translation et dilatation de rapport non nul, et rendant continue l'addition
$$
(x, y) \in (X \times X, \mathcal{F} \otimes \mathcal{F}) \longmapsto x + y \in (X, \mathcal{F})
$$
Remarquons que, comme nous l'avons fait dans le premier complément, il pourrait même suffire que, pour certains ensembles $A \in \mathcal{F}$ (par exemple des convexes),
$$
\big\{ (x, y) \in X \times X \, ; \; x + y \in A \big\} \in \mathcal{F} \otimes \mathcal{F}
$$
Désignons par $\mathcal{C}_0(\mathcal{F})$ l'ensemble de tous les convexes internes mesurables de $X$ contenant $0$ (si $(X, \mathcal{C}_0, \mathcal{F}, \tau)$ est un e.v.l.c.m., $\mathcal{C}_0 \subset \mathcal{C}_0(\mathcal{F})$, mais il n'y a pas égalité en général). Soit $\mathcal{C}_0$ un sous-ensemble de $\mathcal{C}_0(\mathcal{F})$ et $\tau$ la topologie d'espace vectoriel localement convexe engendrée par $\mathcal{C}_0$. Alors, pour toute probabilité $\mu$ sur $\mathcal{F}$, la suite de Cramér associée $(\mu_n)_{n\geqslant 1}$ vérifie un PGD faible relativement à $\tau$. Autrement dit, il y a plus de topologies que celles que nous avons considérées relativement auxquelles on peut énoncer un PGD faible. Si de plus $\mu$ est convexe-tendue, alors, pour toute forme linéaire $\lambda \in X^*$ mesurable, on a
$$
p(\lambda) = \sup_{x \in X} \big( \langle \lambda | x \rangle + s(x) \big)
$$

\subsection{Questions de séparation}

Nous avons supposé l'espace topologique $(X, \tau)$ séparé pour utiliser la notion habituelle de compacité. Toutefois nous avons seulement recours à l'axiome de Borel-Lebesgue vérifié par les compacts, et non à la séparation (cf. lemme \ref{varadcomp}). On dit qu'une partie $Q$ de $X$ est \emph{quasi-compacte} si de tout recouvrement ouvert de $Q$ on peut extraire un sous-recouvrement fini. On dit qu'une partie de $X$ est \emph{relativement quasi-compacte} si elle contenue dans une partie quasi-compacte de $X$. Il est alors très simple d'adapter l'ensemble du texte sans supposer l'espace $(X, \tau)$ séparé : il suffit de remplacer chaque occurrence du mot \og{}compact\fg{} par \og{}quasi-compact\fg{} (dans la borne supérieure faible, dans la définition de la convexe-tension, etc.). Voici ce que devient le lemme de Varadhan convexe. On montre en réalité un résultat plus fort : pour tout $K \in \mathcal{F}$ relativement quasi-compact,
$$
\limsup_{n \to \infty} \frac{1}{n} \log \int_K e^{nf(x)} d\mu_n(x) \leqslant \sup_{Q(K)} [f+s]
$$
où $Q(K)$ est l'intersection de tous les quasi-compacts contenant $K$ (cet ensemble n'est pas nécessairement quasi-compact lui-même). Sachant que l'intersection d'un quasi-compact et d'un fermé est encore quasi-compact, on a $Q(K) \subset \smash{\overline{K}}$ : d'où le résultat et, si $X$ n'est pas séparé, la borne supérieure obtenue peut être bien meilleure.

\subsection{Tension et convexe-tension}

Soient $X_1$ (resp. $X_2$) un espace vectoriel réel, $\mathcal{F}_1$ (resp. $\mathcal{F}_2$) une tribu sur $X_1$ (resp. $X_2$), et $\tau_1$ (resp. $\tau_2$) une topologie séparée sur $X_1$ (resp. $X_2$). On munit $X_1 \times X_2$ de $\mathcal{F}_1 \otimes \mathcal{F}_2$ et de $\tau_1 \times \tau_2$ et on se donne $\mu$ une probabilité sur $X_1 \times X_2$. Si la marginale $\mu_1$ (resp. $\mu_2$) est convexe-tendue sur $D_1$ (resp. $D_2$), alors $\mu$ est convexe-tendue sur $D_1 \times D_2$. En effet, si $(K_{m, 1})_{m \geqslant 1}$ (resp. $(K_{m, 2})_{m \geqslant 1}$) est une suite de parties mesurables de $X_1$ (resp. $X_2$), convexes sur $D_1$ (resp. $D_2$) et relativement compacts sur $D_1$ (resp. $D_2$) telle que $\mu_1(K_{m, 1}) \to 1$ (resp. $\mu_2(K_{m, 2}) \to 1$), alors, pour tout $m \geqslant 1$, $K_{m, 1} \times K_{m, 2}$ est mesurable, convexe sur $D_1 \times D_2$ et relativement compact sur $D_1 \times D_2$, et on a :
$$
\mu(X_1 \times X_2 \setminus K_{m, 1} \times K_{m, 2}) \leqslant \mu_1(X_1 \setminus K_{m, 1}) + \mu_2(X_2 \setminus K_{m, 2}) \to 0
$$

\textbf{Remarque :} Soient $(X, \tilde{\tau})$ un e.v.l.c. quasi-complet et $\tau$ une topologie d'e.v.l.c. compatible avec la dualité entre $X$ et $X^*$ (en fait, $\tau$ est aussi quasi-complète ; cf. \cite[IV.5]{BouEVT}) ; on suppose que $\mathcal{F}$ est la tribu borélienne de $\tau$. Alors, sur $(X, \mathcal{F}, \tau)$, la convexe-tension équivaut à la \emph{tension} : on dit qu'une probabilité $\mu$ est tendue s'il existe une suite $(K_m)_{m \geqslant 1}$ de compacts telle que
$$
\lim_{m\to\infty}\mu(K_m) = 1
$$
Le résultat est conséquence du théorème de Krein (cf. \cite[IV.37, théorème 3]{BouEVT}). Le choix de $\mathcal{F}$ borélienne évite le problème de mesurabilité de l'enveloppe fermée convexe.


\bibliographystyle{alpha-fr}
\bibliography{../cramer}

\end{document}